\documentclass{amsart}

\usepackage{amsthm}
\usepackage{amsfonts}
\usepackage{amsmath}
\usepackage{amssymb}
\usepackage{mathrsfs}
\usepackage{fullpage}
\usepackage{stmaryrd}
\usepackage{hyperref}
\usepackage{verbatim}

\theoremstyle{plain}

\theoremstyle{definition}

\theoremstyle{remark}

\newcommand{\Begin}[2]{\begin{#1}\label{#2}}

\newcommand{\R}{\mathbb{R}}
\newcommand{\Q}{\mathbb{Q}}
\newcommand{\Z}{\mathbb{Z}}

\newcommand{\bPi}{\mathbf{\Pi}}
\newcommand{\bSigma}{\mathbf{\Sigma}}
\newcommand{\bDelta}{\mathbf{\Delta}}
\newcommand{\bbC}{\mathbb{C}}

\newcommand{\bbP}{\mathbb{P}}
\newcommand{\bbQ}{\mathbb{Q}}
\newcommand{\bbR}{\mathbb{R}}

\newcommand{\CM}{\mathcal{M}}
\newcommand{\CN}{\mathcal{N}}

\newcommand{\SCRL}{\mathscr{L}}

\newcommand{\forces}{\Vdash}
\newcommand{\analytic}{{\bSigma_1^1}}
\newcommand{\lanalytic}{{\Sigma_1^1}}
\newcommand{\coanalytic}{{\bPi_1^1}}
\newcommand{\lcoanalytic}{{\Pi_1^1}}
\newcommand{\borel}{{\bDelta_1^1}}
\newcommand{\lborel}{{\Delta_1^1}}
\newcommand{\cantorspace}{{{}^\omega 2}}
\newcommand{\bairespace}{{{}^\omega\omega}}
\newcommand{\finBinarySequence}{{{}^{<\omega}2}}

\newcommand{\F}{{F_{\omega_1}}}
\newcommand{\EBorel}{{E_{\bDelta_1^1}}}

\begin{document}

\title{Equivalence Relations Which Are Borel Somewhere}

\author{William Chan}
\address{Department of Mathematics, California Institute of Technology, Pasadena, CA 91106}
\email{wcchan@caltech.edu}

\begin{abstract}
The following will be shown: Let $I$ be a $\sigma$-ideal on a Polish space $X$ so that the associated forcing of $I^+$ $\borel$ sets ordered by $\subseteq$ is a proper forcing. Let $E$ be an $\analytic$ or a $\coanalytic$ equivalence relation on $X$ with all equivalence classes $\borel$. If for all $z \in H_{(2^{\aleph_0})^+}$, $z^\sharp$ exists, then there exists a $\borel$ set $C \subseteq X$ such that $E \upharpoonright C$ is a $\borel$ equivalence relation. 
\end{abstract}

\maketitle\let\thefootnote\relax\footnote{December 8, 2015
\\*\indent Research partially supported by NSF grants DMS-1464475 and EMSW21-RTG DMS-1044448}


\section{Introduction}\label{Introduction}
The basic question addressed here in its most naive form is: 
\\*
\\*\noindent \textbf{Question:} If $E$ is an equivalence relation on a Polish space $X$, is there a set $C \subseteq X$ such that $E \upharpoonright C$ is a $\bDelta_1^1$ equivalence relation?
\\*
\\*\indent Here, $E \upharpoonright C = E \cap (C \times C)$. It is the substructure of $E$ induced by $C$ when $E$ is considered a structure in the language with a single binary relation symbol.

There are two immediate concerns about how the question is phrased: 

The basic idea of the question is that given an equivalence relation $E$, can one find a subset $C$ such that $E \upharpoonright C$ is a simpler equivalence relation, in particular $\bDelta_1^1$. One does not want to hide any complexity of $E \upharpoonright C$ inside the set $C$. So, the question should be amended to stipulate that $C$ is a $\bDelta_1^1$ subset of $X$. 

Every equivalence relation restricted to a countable set is a $\bDelta_1^1$ equivalence relation. Conditions must be imposed on $C$ to make the question meaningful. $\sigma$-ideals on the Polish space $X$ would include all countable subsets of $X$. So if one demands that $C$ be a non-small set according to a $\sigma$-ideal $I$ on $X$, then the most egregious trivialities vanish. Subsets $C \subseteq X$ with $C \notin I$ are called $I^+$ sets. In the question, a reasonable requirement on $C$ should be that it is $I^+$ and $\borel$.

Without $I$ having some useful properties, there seems to be no particular reason to expect any interesting answer. Some conditions should be imposed on $I$: Given a $\sigma$-ideal $I$ on a Polish space $X$, there is a natural forcing $\bbP_I$ associated with $I$ that has been used extensively in Descriptive Set Theory and Cardinal Characteristics of the Continuum. $\bbP_I$ consists of all $I^+$ $\borel$ subsets of $X$ ordered by $\subseteq$. Motivated by works in Cardinal Characteristics, one could require $I$ to have the property that $\bbP_I$ is a proper forcing. 

In Cardinal Characteristics, properness is used for preservation of certain properties under countable support iterations. This will not be how properness is used in this paper. Rather, properness will be used to produce $I^+$ $\borel$ subsets for which forcing and absoluteness can be used to derive meaningful information. The main tool that makes this approach possible is the following result:
\\*
\\*\noindent\textbf{Fact \ref{properness equivalence}.}
(Zapletal, \cite{Forcing-Idealized} Proposition 2.2.2.) \textit{Let $I$ be a $\sigma$-ideal on a Polish space $X$. The following are equivalent:}

\textit{(i) $\bbP_I$ is a proper forcing.}

\textit{(ii) For any sufficiently large cardinal $\Theta$, for every $B \in \bbP_I$, and for every countable $M \prec H_\Theta$ with $\bbP_I \in M$ and $B \in M$, the set $C := \{x \in B : x \text{ is $\bbP_I$-generic over $M$}\}$ is $I^+$ $\bDelta_1^1$.}
\\*
\\*\indent With this result, the question is now asked with respect to a $\sigma$-ideal such that $\bbP_I$ is proper. 

A natural place to begin exploring this question is with the simplest class of definable equivalence relations just beyond $\borel$ equivalence relations: If $I$ is a $\sigma$-ideal such that $\bbP_I$ is proper and $E$ is an $\analytic$ equivalence relation, is there an $I^+$ $\borel$ set $C$ such that $E \upharpoonright C$ is $\borel$?

Unfortunately, the answer is no.

\Begin{fact}{not analytic arrow borel}
\emph{(\cite{Canonical-Ramsey-Theory-on-Polish-Spaces})} There exists an $\bSigma_1^1$ equivalence relation $E$ and a $\sigma$-ideal $I$ with $\bbP_I$ proper such that for all $I^+$ $\borel$ set $C$, $E \upharpoonright C$ is $\analytic$. 
\end{fact}

\begin{proof}
See \cite{Canonical-Ramsey-Theory-on-Polish-Spaces}, Example 4.25.
\end{proof}

This suggests that in order to possibly obtain a positive answer, the equivalence relation considered should more closely resemble $\borel$ equivalence relations. An obvious feature of $\borel$ equivalence relations is that all their equivalence classes are $\borel$. Kanovei, Sabok, and Zapletal then asked the following question of $\analytic$ equivalence relations which share this feature:

\Begin{question}{introduction main question analytic}
(\cite{Canonical-Ramsey-Theory-on-Polish-Spaces} Question 4.28) If $I$ is a $\sigma$-ideal on a Polish space $X$ such that $\bbP_I$ is proper and $E$ is a $\analytic$ equivalence relation with all classes $\borel$, then is there an $I^+$ $\borel$ set $C$ such that $E \upharpoonright C$ is $\borel$?
\end{question}

Similarly, the question can be asked for the dual class of equivalence relations on the same projective level:

\Begin{question}{introduction main question coanalytic}
If $I$ is a $\sigma$-ideal on a Polish space $X$ such that $\bbP_I$ is proper and $E$ is a $\coanalytic$ equivalence relation with all classes $\borel$, then is there an $I^+$ $\borel$ set $C$ such that $E \upharpoonright C$ is $\borel$?
\end{question}

With these restrictions, the initial naive question becomes a rather robust question. Throughout the paper, the term ``main question'' will refer to questions of the former type for various classes of definable equivalence relations on Polish spaces. For concreteness, the reader should perhaps keep in mind the following explicit instance of the main question: If $E$ is an $\analytic$ equivalence relation with all classes $\borel$, is there a nonmeager or positive measure $\borel$ set $C$ such that $E \upharpoonright C$ is a $\borel$ equivalence relation?

Section \ref{Basic Concepts} will provide the basic concepts from idealized forcing including the main tool about proper idealized forcings used throughout the paper. Some useful notations for expressing the main question is also introduced. The main question in a slightly stronger form is formalized. 

Section \ref{Examples} will provide known results and examples to show that the main question has a positive answer for the most natural $\analytic$ equivalence relations with all $\borel$ classes. In particular, \cite{Canonical-Ramsey-Theory-on-Polish-Spaces} showed that $\analytic$ equivalence relations with all classes countable and equivalence relations $\borel$ reducible to orbit equivalence relations of Polish group actions have a positive answer to the main question. 

Section \ref{Consistency of Positive Answer} will show that a positive answer to the main question for $\analytic$ equivalence relations with all $\borel$ classes follows from some large cardinal assumptions. In particular, it can be proved from iterability principles (such as the existence of a measurable cardinal):
\\*
\\*\noindent\textbf{Theorem \ref{sharps of reals and positive answer main question}.}
\textit{Let $I$ be a $\sigma$-ideal on a Polish space $X$ such that $\bbP_I$ is proper. Let $E$ be an $\analytic$ equivalence relation on $X$ with all classes $\borel$. If for all $z \in \bairespace$, $z^\sharp$ exists and $(\chi_E^I)^\sharp$ exists, then there is an $I^+$ $\borel$ set $C$ such that $E \upharpoonright C$ is $\borel$.}
\\*
\\*\indent Here $\chi_E^I$ is a set depending on $I$ and $E$. This set $\chi_I^E$ is in $H_{(2^{\aleph_0})^+}$ so it is a fairly small set. More explicitly, $\chi_E^I$ is a triple $\langle \bbP_I, \mu_E^{I}, \sigma_E^I\rangle$, where $\mu_E^I, \sigma_E^I \in V^{\bbP_I}$ are names that witness two existential formulas. In fact, these two names can be chosen a bit more constructively using the fullness or maximality property of forcing. In particular, there is a positive answer to the main question for $\analytic$ equivalence relations with all $\borel$ classes if there exists a Ramsey cardinal. 

After showing the positive answer follows from certain large cardinal principles, a natural question would be whether it is consistent relative to some large cardinals that there is a negative answer to the main question. The next sections give partial results for a positive answer using different and weaker consistency assumptions for a restricted class of equivalence relations or ideals. Although these results are inherently interesting, these sections should be understood as an attempt to find situations that can not be used to produce a counterexample to a positive answer to the main question. These results seem to enforce the intuition that a universe with very weak large cardinals may be the ideal place to search for such a counterexample. 

In Section \ref{bsigma Equivalence Relations and Some Ideals}, it will be shown that $I_\text{countable}$, the ideal of countable sets, and $I_{E_0}$, the $\sigma$-ideal generated by $\borel$ sets on which $E_0$ is smooth, will always give a positive answer to the main question for $\analytic$ equivalence relation with all $\borel$ classes. The associated forcings for these two $\sigma$-ideals are Sacks forcing and Prikry-Silver forcing, respectively. The meager ideal, $I_\text{meager}$, and the Lebesgue null ideal, $I_\text{null}$, have associated forcing Cohen forcing and Random forcing. It will be shown that
\\*
\\*\noindent\textbf{Theorem \ref{consistency ZFC consistency null meager arrow borel}.}
\textit{\emph{($\mathsf{ZFC + MA + \neg CH}$)} Let $I$ be either $I_\text{null}$ or $I_\text{meager}$. Let $E$ is an $\analytic$ equivalence relation with all classes $\borel$. Then there exists an $I^+$ $\borel$ set $C$ such that $E \upharpoonright C$ is $\borel$.}
\\*
\\*\indent Section \ref{Thin bSigma11 Equivalence Relations} will consider thin $\analytic$ equivalence relations, i.e., equivalence relations with no perfect set of inequivalent elements. Burgess showed that such equivalence relations have $\aleph_0$ or $\aleph_1$ many equivalence classes. This suggests that the main question for thin $\analytic$ equivalence relations with all $\borel$ classes should be approached combinatorially using covering numbers and the properness of $\bbP_I$. For example assuming $\mathsf{PFA}$, there is a positive answer for all $\sigma$-ideals $I$ with $\bbP_I$ proper and $E$ a thin $\analytic$ equivalence relation with all classes $\borel$. However, the combinatorial approach is not the right one. Using definability ideas, the main question for thin $\analytic$ equivalence relations (even without all $\borel$ classes) has a strong positive answer: 
\\*
\\*\noindent\textbf{Theorem \ref{thin analytic canonicalize ev}.} 
\textit{\emph{($\mathsf{ZFC}$)} If $I$ is a $\sigma$-ideal such that $\bbP_I$ is proper and $E$ is a thin $\analytic$ equivalence relation, then there exists a $I^+$ $\borel$ set $C$ such that $C$ is contained in a single $E$-class.}
\\*
\\*\indent Section \ref{Consistency of Positive Answer for bPi11 Equivalence Relations} will show that a positive answer for $\coanalytic$ equivalence relations with all $\borel$ classes follows from sharps in much the same way as in the $\analytic$ case:
\\*
\\*\textbf{Theorem \ref{sharps of reals and positive answer main question coanalytic}.}
\textit{Let $I$ be a $\sigma$-ideal on a Polish space $X$ such that $\bbP_I$ is proper. Let $E$ be a $\coanalytic$ equivalence relation on $X$ with all classes $\borel$. If for all $z \in \bairespace$, $z^\sharp$ exists and $(\chi_E^I)^\sharp$ exists, then there is a $I^+$ $\borel$ set $C$ such that $E \upharpoonright C$ is $\borel$.}
\\*
\\*\indent The set $\chi_E^I$ is defined similarly to the $\analytic$ case. 

Section \ref{bPi11 Equivalence Relations with Thin or Countable Classes} will consider $\coanalytic$ equivalence relations with all classes countable. As mentioned above, $\mathsf{ZFC}$ can provide a positive answer to the main question for $\analytic$ equivalence relation with all countable classes. In the $\bPi_1^1$ case, there is insufficient absoluteness to carry out the same proof. However, from the consistency of a remarkable cardinal, one can obtain the consistency of a positive answer to the main question for $\coanalytic$ equivalence relation with all countable classes:
\\*
\\*\noindent\textbf{Theorem \ref{remarkable cardinal consistency pi11aleph0 arrow borel}.}
\textit{Let $\kappa$ be a remarkable cardinal in $L$. Let $G \subseteq \text{Coll}(\omega, < \kappa)$ be $\text{Coll}(\omega, < \kappa)$ generic over $L$. In $L[G]$, if $I$ is $\sigma$-ideal with $\bbP_I$ proper and $E$ is a $\coanalytic$ equivalence relation with all classes countable, then there exists some $I^+$ $\borel$ set $C$ such that $E \upharpoonright C$ is $\borel$.}
\\*
\\*\indent There is also a similar result using a weakly compact cardinal but $\bbP_I$ must be a $\aleph_1$-c.c. forcing. 

Section \ref{bPi21 Equivalence Relations with all bDelta11 Classes} will show that in $L$, the main question for $\bDelta_2^1$ equivalence relations with all classes $\borel$ (in fact countable) is false. 
\\*
\\*\noindent\textbf{Theorem \ref{not EL canonicalize borel}.}
\textit{In $L$, there is a $\bDelta_2^1$ equivalence relation with all classes countable such that for all $\sigma$-ideals $I$ and all $I^+$ $\borel$ sets $C$, $E \upharpoonright C$ is not $\borel$.}
\\*
\\*\indent It is not known whether a positive answer in this case is consistent; however, there seems to be no reason it could be. 

Finally, the last section will summarize the work of the paper from the point of view of showing the consistency of a negative answer to the main question. Related questions will be introduced. Some dubious speculations about how a negative answer could be obtained will be discussed. 

Drucker, in \cite{Borel-Canonization-of-Analytic-Sets-with-Borel-Sections}, has independently obtained some results that are very similar to what appears in this paper: He has shown that a positive answer to the main question follows from a measurable cardinal using similar ideas to those appearing in Section \ref{Consistency of Positive Answer}. He has obtained results for $\sigma$-ideals whose forcings are provably $\aleph_1$-c.c. which are similar to Section \ref{bsigma Equivalence Relations and Some Ideals}. Drucker also proved the results of Section \ref{bPi21 Equivalence Relations with all bDelta11 Classes} of this paper using a very similar equivalence relation. In \cite{Borel-Canonization-of-Analytic-Sets-with-Borel-Sections}, Drucker also considers more general forms of canonicalization than what appears in this paper.

The author would like to thank Ohad Drucker, Sy-David Friedman, and Alexander Kechris for many helpful discussions about this paper.

\section{Basic Concepts}\label{Basic Concepts}
This section reviews the basics of idealized forcing and formalizes the main question of interest. 

\Begin{definition}{forcing from ideal}
Let $I$ be a $\sigma$-ideal on a Polish space $X$. Let $\bbP_I$ be the collection of all $I^+$ $\bDelta_1^1$ subsets of $X$. Let $\leq_{\bbP_I} = \subseteq$. Let $1_{\bbP_I} = X$. $(\bbP_I, \leq_{\bbP_I}, 1_{\bbP_I})$ is the forcing associated with the ideal $I$. 
\end{definition}

\Begin{fact}{canonical name generic point}
Let $I$ be a $\sigma$-ideal on a Polish space $X$. There is a $\bbP_I$-name $\dot x_\textup{gen}$ such that for all $\bbP_I$-generic filters $G$ over $V$ and all $B$ which is $\borel$ coded in $V$, $B \in G$ if and only if $\dot x_\textup{gen}[G] \in B$. 
\end{fact}

\begin{proof}
See \cite{Forcing-Idealized}, Proposition 2.1.2.
\end{proof}

\Begin{definition}{generic point of ideal forcing}
Let $I$ be a $\sigma$-ideal on a Polish space $X$. Let $M \prec H_\Theta$ be a countable elementary substructure for some cardinal $\Theta$. $x \in X$ is $\bbP_I$-generic over $M$ if and only if the set $\{A \in \bbP_I \cap M : x \in A\}$ is a $\bbP_I$-generic filter over $M$. 
\end{definition}

\Begin{fact}{properness equivalence}
Let $I$ be a $\sigma$-ideal on a Polish space $X$. The following are equivalent: 

(i) $\bbP_I$ is a proper forcing.

(ii) For any sufficiently large cardinal $\Theta$, for every $B \in \bbP_I$, and for every countable $M \prec H_\Theta$ with $\bbP_I \in M$ and $B \in M$, the set $C := \{x \in B : x \text{ is $\bbP_I$-generic over $M$}\}$ is $I^+$ $\bDelta_1^1$.
\end{fact}

\begin{proof}
See \cite{Forcing-Idealized}, Proposition 2.2.2. Since this is the most important tool in this paper, a proof will be sketched: 

(i) $\Rightarrow$ (ii) Let $B \in \bbP_I \cap M$ be arbitrary. It is straightforward to show that $C$ is $\borel$. Suppose $C \in I$. Then by Fact \ref{canonical name generic point}, $B \forces_{\bbP_I}^V \dot x_\text{gen} \notin C$. This implies that there is some $D \subseteq \bbP_I$ which is dense, $D \in M$, and $B \forces^V_{\bbP_I} \check M \cap \check D \cap \dot G = \emptyset$. Therefore, there can be no $(M,\bbP_I)$-generic condition below $B$. $\bbP_I$ is not proper. 

(ii) $\Rightarrow$ (i) Let $B \in \bbP_I \cap B$ be arbitrary. Suppose $C \notin I$. Then $C \forces_{\bbP_I}^V \dot x_\text{gen} \in C$. So for all $D \subseteq \bbP_I$ with $D$ dense and $D \in M$, $C \forces_{\bbP_I}^V \check M \cap \check D \cap \dot G \neq \emptyset$. $C$ is a $(M,\bbP_I)$-generic condition below $B$. $\bbP_I$ is proper.
\end{proof}

The following is some convenient notation:

\Begin{definition}{arrow notation}
(\cite{Canonical-Ramsey-Theory-on-Polish-Spaces} Definition 1.15) Let $\Lambda$ and $\Gamma$ be classes of equivalence relations defined on $\bDelta_1^1$ subsets of Polish spaces. Let $I$ be a $\sigma$-ideal on a Polish space $X$. Define $\Lambda \rightarrow_I \Gamma$ to mean: for all $B$ which are $I^+$ $\bDelta_1^1$ subsets of $X$ and all equivalence relation $E$ defined on $X$ such that $E \upharpoonright B \in \Lambda$, there exists a $I^+$ $\borel$ set $C \subseteq B$ such that $E \upharpoonright C \in \Gamma$. 
\end{definition}

The following are some of the classes of equivalence relations that will appear later.

\Begin{definition}{some classes of equivalence relation}
For any Polish space $X$, $\text{ev}$ denote the full equivalence relation on $X$ consisting of a single class.

For any Polish space $X$, $\text{id}$ is the equality equivalence relation. 

$\borel$ denote the class of all $\borel$ equivalence relations defined on $\borel$ subsets of Polish spaces. (In context, it should be clear when $\borel$ refers to the class of equivalence relations or just the $\borel$ definable subsets.)

$\analytic^\borel$ is the class of all $\analytic$ equivalence relations defined on $\borel$ subsets of Polish spaces with all classes $\borel$.

$\coanalytic^\borel$ is the class of all $\coanalytic$ equivalence relations defined on $\borel$ subsets of Polish spaces with all classes $\borel$.

$\bDelta_2^{1 \ \borel}$ is the class of all $\bDelta_2^1$ equivalence relations defined on $\borel$ subsets of Polish spaces with all classes $\borel$. 

A thin equivalence relation is an equivalence relation with no perfect set of inequivalent elements. 

$\bSigma_1^{1 \ \text{thin}}$ is the class of all thin $\analytic$ equivalence relations defined on $\borel$ subsets of Polish spaces. 

$\bSigma_1^{1 \ \text{thin} \ \borel}$ is the class of all thin $\analytic$ equivalence relations defined on $\borel$ subsets of Polish spaces and have all classes $\borel$. 

$\bPi_1^{1 \ \aleph_0}$ denotes the class of all $\coanalytic$ equivalence relations with all classes countable and defined on $\borel$ subsets of Polish spaces. 

A thin set is a set without a perfect subset. 

Let $\bPi_1^{1 \ \text{thin}}$ denote the class of all $\coanalytic$ equivalence relations with all classes thin and defined on $\borel$ subsets of Polish spaces. 
\end{definition}

Kanovei, Sabok, and Zapletal asked the following questions:

\Begin{question}{main question}
(\cite{Canonical-Ramsey-Theory-on-Polish-Spaces} Question 4.28) If $I$ is a $\sigma$-ideal on a Polish space $X$ such that $\bbP_I$ is proper, then does $\analytic^\borel \rightarrow_I \bDelta_1^1$ hold?
\end{question}

This paper will address this question and its various related forms for other classes of definable equivalence relations.

\section{Examples}\label{Examples}

This section gives known results concerning the main question and some examples.

\Begin{proposition}{reducible orbit equiv polish group action arrow borel}
Let $\Gamma_1$ denote the class of equivalence relation $\bDelta_1^1$ reducible to orbit equivalence relation of Polish group actions. Then $\Gamma_1 \rightarrow_I \bDelta_1^1$ for any $\sigma$-ideal $I$ on $X$ such that $\bbP_I$ is proper.  
\end{proposition}

\begin{proof}
See \cite{Canonical-Ramsey-Theory-on-Polish-Spaces}, Theorem 4.26. 
\end{proof}

\Begin{proposition}{all classes countable arrow borel}
Let $\Gamma_2$ denote the class of $\bSigma_1^1$ equivalence relation with all classes countable. Then $\Gamma_2 \rightarrow_I \bDelta_1^1$ for any $\sigma$-ideal $I$ on $X$ such that $\bbP_I$ is proper. 
\end{proposition}

\begin{proof}
See \cite{Canonical-Ramsey-Theory-on-Polish-Spaces}, Theorem 4.27. The proof is provided below to emphasize a particular observation.

Fix a $B \subseteq X$ which is $I^+$ $\bDelta_1^1$. As $E$ is $\bSigma_1^1$, there exists some $z \in \cantorspace$ such that $E$ is $\Sigma_1^1(z)$. For each $x \in X$, $[x]_E$ is $\Sigma_1^1(x,z)$. Since every $\Sigma_1^1(x,z)$ set with a non-$\Delta_1^1(x,z)$ element has a perfect subset (see \cite{Recursive-Aspects-of-Descriptive-Set-Theory} Theorem 6.3), the statement ``all $E$-classes are countable'' is equivalent to 
$$(\forall x)(\forall y)(y \ E \ x \Rightarrow y \in \Delta_1^1(x,z))$$
As the relation ``$y \in \Delta_1^1(x,z)$'' in variable $x$ and $y$ is $\Pi_1^1(z)$, the above is $\Pi_1^1(z)$. By Mostowski absoluteness, $1_{\bbP_I} \forces_{\bbP_I}$ ``All $E$-classes are countable''. There is some $\bbP_I$-name $\tau$ such that $B \forces_{\bbP_I} \tau \in {}^\omega X \wedge \tau \text{ enumerates } [\dot x_\text{gen}]_E$. By \cite{Forcing-Idealized} Proposition 2.3.1, there exists some $B' \subseteq B$ with $B' \in \bbP_I$ and a $\bDelta_1^1$ function $f$ such that $B' \forces_{\bbP_I} f(\dot x_\text{gen}) = \tau$. Choose $M \prec H_\Theta$ with $\Theta$ sufficiently large and $M$ contains $\bbP_I$, $B'$, $\tau$, and the code for $f$. By Fact \ref{properness equivalence}, let $C \subseteq B'$ be the $I^+$ $\bDelta_1^1$ set of $\bbP_I$-generic over $M$ elements in $B'$. 

The claim is that for $x,y \in C$, $x \ E \ y$ if and only if $(\exists n)(f(n) = y)$. This is because for all $x \in C$, $M[x] \models f(x) \text{ enumerates } [x]_E$. Let $N$ be the Mostowski collapse of $M[x]$. One can always assume the transitive closure of elements of $X$ is a subset of $M$ (for instance, one could have identified $X$ with $\bairespace$). Therefore the Mostowski collapse map does not move elements of $X$. Hence $N \models f(x) \text{ enumerates } [x]_E$. The statement ``$f(x)$ enumerates $[x]_E$'' is the conjunction of a $\bSigma_1^1$ and $\bPi_1^1$ formula coded in $N$. By Mostowski absoluteness, $f(x)$ enumerates $[x]_E$ in $V$. This proves the claim. Thus $E \upharpoonright C$ is a $\bDelta_1^1$ equivalence relation. 
\end{proof}

By the work above, for each $x$ which is $\bbP_I$-generic over $M$, $M[x] \models [x]_E \text{ is countable}$. So in $M[x]$, there exists some real $u_x$ such that $u_x$ codes an enumeration of $[x]_E$. In $M[x]$, $[x]_E$ is $\Delta_1^1(u_x)$. In the above proof, one showed that $u_x$ remains an enumeration of $[x]_E$ even in $V$. So $[x]_E$ is $\Delta_1^1(u_x)$ even in $V$. This observation is the quintessential idea of the proof of the positive answer for the main question assuming large cardinal properties. Note that in the above proof, there was a $\bDelta_1^1$ function $f$ which uniformly provided the enumeration of $[x]_E$ for each $x \in C$. This feature is not necessary.

Below, positive answers to the main question will be demonstrated for some specific equivalence relations.

\Begin{definition}{admissible ordinal equivalence relation}
For each $x \in \cantorspace$, $\omega_1^x$ is the least $x$-admissible ordinal above $\omega$. Define the equivalence relation $\F$ on $\cantorspace$ by $x \ \F \ y$ if and only if $\omega_1^x = \omega_1^y$. 

$\F$ is an $\Sigma_1^1$ equivalence relation with all classes $\bDelta_1^1$. 
\end{definition}

\Begin{example}{F equivalence relation arrow borel}
Let $I$ be a $\sigma$-ideal on $\cantorspace$ with $\bbP_I$ proper. Then $\{\F\} \rightarrow_I \{\text{ev}\}$, i.e. there is an $I^+$ class. 
\end{example}

\begin{proof}
Let $B$ be an arbitrary $I^+$ $\bDelta_1^1$ set. Choose $M \prec H_\Theta$ where $\Theta$ is a sufficiently large cardinal and $\bbP_I, B \in M$. By Fact \ref{properness equivalence}, let $C \subseteq B$ be the $I^+$ $\bDelta_1^1$ set of $\bbP_I$-generic over $M$ reals in $B$. For each $x \in C$, $\omega_1^x \in M[x] \cap \text{ON}$. Since the ground model and the forcing extension have the same ordinals, $\omega_1^x \in M \cap \text{ON}$. For each $\alpha \in \text{ON}$, let $F_{\omega_1}^\alpha = \{x \in \cantorspace : \omega_1^x = \alpha\}$. Each $F_{\omega_1}^\alpha = \emptyset$ or is an $\F$-class. $C = \bigcup_{\alpha \in M \cap \text{ON}} F_{\omega_1}^\alpha \cap C$. Since $\F$-classes are $\bDelta_1^1$, $F_{\omega_1}^\alpha \cap C$ is $\bDelta_1^1$ for all $\alpha$. As $M$ is countable, $M \cap \text{ON}$ is countable. There exists some $\alpha \in M \cap \text{ON}$ such that $F_{\omega_1}^\alpha \cap C$ is $I^+$ since $I$ is a $\sigma$-ideal. For this $\alpha$, $\F \upharpoonright F_{\omega_1}^\alpha \cap C = \text{ev} \upharpoonright F_{\omega_1}^\alpha \cap C$. 
\end{proof}

Since $M \prec H_\Theta$, for each $x \in C$, there exists a countable admissible ordinal $\alpha > \omega_1^x$ with $\alpha \in M$. By Sacks' theorem applied in $M$, let $y \in M$ be such that $\omega_1^y = \alpha$. Then $[x]_\F$ is $\Delta_1^1(y)$. Again the phenomenon described above occurs: there exist some $y \in M[x]$ (in fact $y \in M$) such that $M[x] \models [x]_\F \text{ is } \Delta_1^1(y)$ and $V \models [x]_\F \text{ is } \Delta_1^1(y)$. 

Actually, $\F$ is classifiable by countable structures. Proposition \ref{reducible orbit equiv polish group action arrow borel} would have already shown $\{\F\} \rightarrow_I \bDelta_1^1$. See \cite{The-Countable-Admissible-Ordinal-Equivalence-Relation} for more information about $\F$. 

\Begin{definition}{Ew1 equivalence relation}
Define the equivalence relation $E_{\omega_1}$ on $\cantorspace$ by
$$x \ E_{\omega_1} \ y \Leftrightarrow (x \notin \text{WO} \wedge y \notin \text{WO}) \vee (\text{ot}(x) = \text{ot}(y))$$
where $\text{WO}$ is the set of reals coding well-orderings and for $x \in \text{WO}$, $\text{ot}(x)$ is the order type of the linear order coded by $x$.

$E_{\omega_1}$ is a $\Sigma_1^1$ equivalence relation with all classes $\bDelta_1^1$ except for one $\Sigma_1^1$ class consisting of the reals that do not code wellfounded linear orderings.
\end{definition}

\Begin{example}{Ew1 arrow borel}
Let $I$ be a $\sigma$-ideal on $\cantorspace$ with $\bbP_I$ proper. Then $\{E_{\omega_1}\} \rightarrow_I \{\text{ev}\}$. 
\end{example}

\begin{proof}
Let $B \subseteq \cantorspace$ be $I^+$ $\bDelta_1^1$. 

(Case I) There exists some $B' \leq_{\bbP_I} B$ such that $B' \forces_{\bbP_I} \dot x_\text{gen} \notin \text{WO}$: Let $M \prec H_\Theta$ be a countable elementary structure with $\Theta$ a sufficiently large cardinal and $\bbP_I, B' \in M$. By Fact \ref{properness equivalence}, let $C \subseteq B'$ be the $I^+$ $\bDelta_1^1$ set of all $\bbP_I$-generic over $M$ reals in $B'$. Let $x \in C$. By Fact \ref{canonical name generic point}, let $G_x \subseteq \bbP_I$ be the generic filter associated with $x$. $B' \in G_x$ since $x \in B'$. $B' \forces_{\bbP_I} \dot x_\text{gen} \notin \text{WO}$ implies that $M[x] \models x \notin \text{WO}$. Let $N$ be the Mostowski collapse of $M[x]$. Since the Mostowski collapse map does not move reals, $x \in N$. Also $N \models x \notin \text{WO}$. Since $\text{WO}$ is $\Pi_1^1$, $V \models x \notin \text{WO}$. Hence $E \upharpoonright C$ consists of a single class. So $E \upharpoonright C = \text{ev} \upharpoonright C$. 

(Case II) $B \forces_{\bbP_I} \dot x_\text{gen} \in \text{WO}$: Then choose $M \prec H_\Theta$ a countable elementary substructure and $\Theta$ a sufficiently large cardinal. By Fact \ref{properness equivalence}, let $C \subseteq B$ be the $I^+$ $\bDelta_1^1$ set of $\bbP_I$-generic over $M$ reals in $B$. As in Case I, $B \forces_{\bbP_I} \dot x_\text{gen} \in \text{WO}$ implies that $V \models x \in \text{WO}$. So when $x \in C$, $\text{ot}(x) \in M[x] \cap \text{ON}$. For each ordinal $\alpha < \omega_1$, $E_{\omega_1}^\alpha := \{x \in \text{WO} : \text{ot}(x) = \alpha\}$ is a $\bDelta_1^1$ set. Since $M[x]$ and $M$ have the same ordinals and $M$ is countable, $M[x]$ has only countably many ordinals. $C = \bigcup_{\alpha \in M \cap \text{ON}} E_{\omega_1}^\alpha \cap C$. $E_{\omega_1}^\alpha \cap C$ is $\bDelta_1^1$ for each $\alpha$. Since $I$ is a $\sigma$-ideal, there is some $\alpha \in M \cap \text{ON}$ such that $E_{\omega_1}^\alpha \cap C$ is $I^+$. So for this $\alpha$, $E_{\omega_1} \upharpoonright E_{\omega_1}^\alpha \cap C = \text{ev} \upharpoonright E_{\omega_1}^\alpha \cap C$. 
\end{proof}

Note that $E_{\omega_1}$ does not have all classes $\bDelta_1^1$. However, it is a thin $\Sigma_1^1$ equivalence relation. It will be shown later that the main question can be answered positively for thin $\bSigma_1^1$ equivalence relation regardless of whether the classes are all $\bDelta_1^1$. 

Next, there is one further enlightening example which does not fall under the scope of Proposition \ref{reducible orbit equiv polish group action arrow borel} or Proposition \ref{all classes countable arrow borel}.

\Begin{fact}{universal for borel equivalence relation}
There exist a $\Pi_1^1$ set $D \subseteq \bairespace$, a $\Pi_1^1$ set $P \subseteq (\bairespace)^3$, and a $\Sigma_1^1$ set $S \subseteq (\bairespace)^3$ such that

(1) If $z \in D$, then for all $x,y \in \bairespace$, $P(z,x,y) \Leftrightarrow S(z,x,y)$.

(2) For all $z \in D$, the relation $x \ E_z \ y$ if and only $P(z,x,y)$ is an equivalence relation, which is $\bDelta_1^1$ equivalence relation by (1).

(3) If $E$ is a $\bDelta_1^1$ equivalence relation, then there is a $z$ such that $x \ E \ y \Leftrightarrow P(z,x,y)$. 
\end{fact}

\begin{proof}
See \cite{Equivalence-Relations-which-Reduce-all-Borel-Equivalence-Relations}, Definition 14.
\end{proof}

\Begin{definition}{equiv containing all borel equiv}
(\cite{Equivalence-Relations-which-Reduce-all-Borel-Equivalence-Relations} Definition 29) Define the equivalence relation $\EBorel$ on $(\bairespace)^2$ by
$$(z_1, x_1) \ \EBorel \ (z_2, x_2) \Leftrightarrow (z_1 = z_2) \wedge (\neg D(z_1) \vee S(z_1, x_1, x_2))$$

$\EBorel$ is a $\Sigma_1^1$ equivalence relation. For each $z$, define $E_z$ by $x \ E_z \ y$ if and only if $(z,x) \ \EBorel \ (z,y)$. For all $z$, $E_z$ is a $\bDelta_1^1$ equivalence relation. If $E$ is a $\bDelta_1^1$ equivalence relation, then there exists a $z \in D$ such that $E = E_z$. $\EBorel$ has all classes $\bDelta_1^1$. 
\end{definition}

\Begin{fact}{borel equivalence relation reduce into Eborel}
If $E$ is a $\bDelta_1^1$ equivalence relation, then $E \leq_{\bDelta_1^1} \EBorel$. 
\end{fact}

\begin{proof}
(See \cite{Equivalence-Relations-which-Reduce-all-Borel-Equivalence-Relations}) Let $z \in D$ such that $x \ E \ y \Leftrightarrow S(z,x,y)$ for all $x,y \in \bairespace$. Then $f : \bairespace \rightarrow (\bairespace)^2$ defined by $f(x) = (z,x)$ is the desired reduction.
\end{proof}

\Begin{proposition}{EBorel is not reducible to orbit equiv}
$\EBorel$ is a $\Sigma_1^1$ equivalence relation with all classes $\bDelta_1^1$, has uncountable classes, and is not reducible to the orbit equivalence relation of a Polish group action.
\end{proposition}

\begin{proof}
All except the last statement have been mentioned above. For the last statement, $E_1$ is a $\Delta_1^1$ equivalence relation so by Fact \ref{borel equivalence relation reduce into Eborel}, $E_1 \leq_{\bDelta_1^1} \EBorel$. $E_1$ is not $\bDelta_1^1$ reducible to any orbit equivalence relation of a Polish group action, by \cite{The-Classification-of-Hypersmooth-Borel-Equivalence-Relations} Theorem 4.2.
\end{proof}

$\EBorel$ is a $\Sigma_1^1$ equivalence relation which does not fall under Proposition \ref{reducible orbit equiv polish group action arrow borel} or Proposition \ref{all classes countable arrow borel}. Next, it will be shown that the main question formulated for $\EBorel$ has a positive answer.

\Begin{theorem}{EBorel arrow borel}
Let $I$ be a $\sigma$-ideal on $(\bairespace)^2$ such that $\bbP_I$ is a proper forcing, then $\{\EBorel\} \rightarrow_I \bDelta_1^1$. 
\end{theorem}

\begin{proof}
Since $D$ is $\Pi_1^1$, let $T$ be a recursive tree on $\omega \times \omega$ such that $x \in D$ if and only $T^x$ is well-founded. Define
$$D_\alpha := \{x : \text{rk}(T^x) < \alpha\}$$
Each $D_\alpha$ is $\bDelta_1^1$. Define
$$(z_1, x_1) \ E_{\bDelta_1^1}^\alpha (z_2, x_2) \Leftrightarrow z_1 = z_2 \wedge (\neg D_\alpha(z_1) \vee P(z_1, x_1, x_2)) \Leftrightarrow z_1 = z_2 \wedge (\neg D_\alpha(z_1) \vee (S(z_1, x_1, x_2)))$$
Each $E_{\bDelta_1^1}^\alpha$ is $\bDelta_1^1$ and $\EBorel = \bigcap_{\alpha < \omega_1} E_{\bDelta_1^1}^\alpha$.

Let $B \subseteq (\bairespace)^2$ be a $\bDelta_1^1$ $I^+$ set.  Let $\pi_1 : (\bairespace)^2 \rightarrow \bairespace$ be the projection onto the first coordinate.

(Case I) $B \not\forces_{\bbP_I} \pi_1(\dot x_\text{gen}) \in D$: Then there exists some $B' \leq_{\bbP_I} B$ such that $B' \forces_{\bbP_I} \pi_1(\dot x_\text{gen}) \notin D$. Now let $M \prec H_\Theta$ be a countable elementary substructure with $B', \bbP_I \in M$ and $\Theta$ some sufficiently large cardinal. By Fact \ref{properness equivalence}, let $C \subseteq B'$ be the $I^+$ $\bDelta_1^1$ set of $\bbP_I$-generic over $M$ elements in $B'$. By elementarity, $B' \forces^M_{\bbP_I} \pi_1(\dot x_\text{gen}) \notin D$. So for all $x \in C$, $M[x] \models \pi_1(x) \notin D$. For each $x \in C$, let $N_x$ denote the Mostowski collapse of $M[x]$. Note that the Mostowski collapse map does not move reals. Hence $N_x \models \pi_1(x) \notin D$. By Mostowski absoluteness, $\pi_1(x) \notin D$. So for all $(z_1, x_1), (z_2, x_2)\in C$, $(z_1,x_2) \ \EBorel (z_2, x_2) \Leftrightarrow z_1 = z_2$. So $\EBorel \upharpoonright C$ is $\bDelta_1^1$. 

(Case II) Otherwise $B \forces_{\bbP_I} \pi_1(\dot x_\text{gen}) \in D$: Let $M \prec H_\Theta$ be a countable elementary substructure with $B, \bbP_I \in M$ and $\Theta$ some sufficiently large cardinal. By Fact \ref{properness equivalence}, let $C \subseteq B$ be the $I^+$ $\bDelta_1^1$ set of $\bbP_I$-generic over $M$ elements in $B$. By elementarity, $B \forces_{\bbP_I}^M \pi_1(\dot x_\text{gen}) \in D$. For all $x \in C$, $M[x] \models \pi_1(x) \in D$. $M[x] \models (\exists \alpha < \omega_1)(\text{rk}(T^{\pi_1(x)}) < \alpha)$. Let $\beta = M \cap \omega_1$. Since $\bbP_I$ preserves $\aleph_1$, $\omega_1^{M[x]} = \omega_1^M$ for all $x \in C$. For each $x \in C$, let $N_x$ be the Mostowski collapse of $M[x]$. Note that the Mostowski collapse map does not move any reals. Then for all $x \in C$, $N_x \models (\exists \alpha < \omega_1^{N_x} = \beta)(\text{rk}(T^{\pi_1(x)}) < \alpha)$. For each $x \in C$, there is some $\alpha < (\omega_1)^{N_x}$ such that $N_x \models \text{rk}(T^{\pi_1(x)}) < \alpha$. After expressing this statement using a real in $N_x$ that code the countable (in $N_x$) ordinal $\alpha$, Mostowski absoluteness implies that $\text{rk}(T^{\pi_1(x)}) < \alpha$. It has been shown that for all $x \in C$, $\text{rk}(T^{\pi_1(x)}) < \beta$. For all $x \in C$, $\pi_1(x) \in D_\beta$. $\EBorel \upharpoonright C = E_{\bDelta_1^1}^\beta \upharpoonright C$. The latter is $\bDelta_1^1$. 
\end{proof}

The above proof motivates the ideas used in the next section.

\section{Positive Answer for $\bSigma_1^{1 \ \borel}$}\label{Consistency of Positive Answer}
Using some of the ideas from the earlier examples, it will be shown that a positive answer to the main question follows from large cardinals. Avoiding any explicit mention of iteration principles, a crude result for the positive answer is first given assuming some generic absoluteness and the existence of tree representations that behave very nicely with generic extensions. This result will illustrate all the main ideas before going into the more optimal but far more technical proof using interable structures.

For simplicity, assume that $E$ is a $\bSigma_1^1$ equivalence relation on $\bairespace$.

First, a classical result about $\bSigma_1^1$ equivalence relations. 

\Begin{fact}{Burgess theorem}
Let $E$ be a $\Sigma_1^1(z)$ equivalence relation on $\bairespace$. Then there exists $\bDelta_1^1$ relations $E_\alpha$, for $\alpha < \omega_1$, with the property that if $\alpha < \beta$, then $E_\alpha \supseteq E_\beta$, $E = \bigcap_{\alpha < \omega_1} E_\alpha$, and there exists a club set $C \subseteq \omega_1$ such that for all $\alpha \in C$, $E_\alpha$ is an equivalence relation.
\end{fact}

\begin{proof}
See \cite{Descriptive-Set-Theory-and-Infinitary-Languages}. Since $E$ is $\Sigma_1^1(z)$, let $T$ be a $z$-recursive tree on $\omega \times \omega \times \omega$ such that $(x,y) \in E$ if and only if $T^{(x,y)}$ is illfounded. For each $\alpha < \omega_1$, define $E_\alpha$ by $(x,y) \in E_\alpha \Leftrightarrow \text{rk}(T^{(x,y)}) > \alpha$. Observe that $E_\alpha$ is $\Delta_1^1(z, c)$ for any $c$ which codes the ordinal $\alpha$. 

The verification of the rest of the theorem is an application of the boundedness theorem and can be found in any reference on the descriptive set theory of equivalence relations. (See Lemma \ref{pi11 increasing union borel equiv} for a similar result in the $\bPi_1^1$ case.)
\end{proof}

For the rest of this section, fix a $z$-recursive tree $T$ as in the proof above. $\{E_\alpha: \alpha < \omega_1\}$ will refer to the sequence of $\bDelta_1^1$ equivalence relations obtained from $T$.

\Begin{lemma}{analytic equivalence stabilization of classes}
Let $E$ be a $\Sigma_1^1(z)$ equivalence relation. Let $x,y \in \bairespace$ be such that $[x]_E$ is $\Pi_1^1(y)$. Let $\delta$ be an ordinal such that $\omega_1^{x \oplus y \oplus z} \leq \delta$ and $E_\delta$ is an equivalence relation. Then $[x]_E = [x]_{E_\delta}$.
\end{lemma}

\begin{proof}
Define $E' \subseteq (\bairespace)^2$ by
$$a \ E' \ b \Leftrightarrow (a \in [x]_E \wedge b \in [x]_E) \vee (a \notin [x]_E \wedge b \notin [x]_E)$$ 
$E'$ is $\Pi_1^1(x \oplus y \oplus z)$. $(\bairespace)^2 - E'$ is then $\Sigma_1^1(x \oplus y \oplus z)$. $(\bairespace)^2 - E' \subseteq (\bairespace)^2 - E$. By the effective boundedness theorem, there exists an $\alpha < \omega_1^{x \oplus y \oplus z} < \delta$ such that for all $(x,y) \in (\bairespace)^2 - E'$, $\text{rk}(T^{(x,y)}) \leq \alpha$. Hence $E' \supseteq E_\alpha$. 

Since $E \subseteq E_\delta$, $[x]_E \subseteq [x]_{E_\delta}$. Since $E_\delta \subseteq E_\alpha \subseteq E'$, $[x]_{E_\delta} = \{u : (u,x) \in E_\delta\} \subseteq \{u : (u,x) \in E'\} = [x]_{E'} = [x]_E$. Therefore, $[x]_E = [x]_{E_\delta}$. 
\end{proof}

Lemma \ref{analytic equivalence stabilization of classes} gives an upper bound on the ordinal level of the sequence $\{E_\alpha : \alpha < \omega_1\}$ where a $\bPi_1^1$ $E$-class stabilizes. According to this lemma, a crucial piece of information in finding this bound is a real $y$ which can be used as a parameter in some $\Pi_1^1$ definition of the $\bPi_1^1$ $E$-class. Rather than knowing a particular $\Pi_1^1$ code, it suffices to know where some particular code lives:

\Begin{lemma}{true pi11 code in generic extensions}
Let $E$ be a $\Sigma_1^1(z)$ equivalence relation. Let $M \prec H_\Theta$ be a countable elementary substructure with $z \in M$. Let $\bbP$ be a forcing in $M$ which adds a generic real. Suppose for all $g$ which are $\bbP$-generic over $M$, there exists a $y \in M[g]$ such that $V \models [g]_E$ is $\Pi_1^1(y)$. Then there exists a countable ordinal $\alpha$ such that for all $\bbP$-generic over $M$ reals $g$, $[g]_E = [g]_{E_\alpha}$. 
\end{lemma}

\begin{proof}
Let $M'$ be the Mostowski collapse of $M$. Let $\alpha = M' \cap \text{ON}$. Since $M$ is countable, $\alpha < \omega_1$. 

Let $g$ be a $\bbP$-generic over $M$ real. Let $y \in M[g]$ be such that $[g]_E$ is $\Pi_1^1(y)$, in $V$. Let $N$ be the Mostowski collapse of $M[g]$. Note that $\alpha = N \cap \text{ON}$. Since $H_\Theta$ is an admissible set, $H_\Theta \models \mathsf{KP}$. $M \models \mathsf{KP}$. $M[g] \models \mathsf{KP}$. So $N \models \mathsf{KP}$. $N$ is an admissible set. $g$, $y$, and $z$ are elements of $N$ since reals are collapsed to themselves. $\alpha$ is a $(g \oplus y \oplus z)$-admissible ordinal. In particular, $\omega_1^{g \oplus y \oplus z} \leq \alpha$.

Next to show that $E_\alpha$ is an equivalence relation: Let $\{\kappa_\eta : \eta < \xi\}$ denote an increasing enumeration of the cardinals of $M'$. If $\eta$ is a limit ordinal and $\beta_\gamma$ has been defined for $\gamma < \eta$ in such a way that $\kappa_\gamma < \beta_\gamma < \kappa_{\gamma + 1}$ and $E_{\beta_\gamma}$ is an equivalence relation (in $V$), then let $\beta_\eta = \sup_{\gamma < \eta} \beta_\gamma$. $E_{\beta_\eta} = \bigcap_{\gamma < \eta} E_{\beta_\gamma}$ so $E_{\beta_\eta}$ is an equivalence relation. Suppose $\eta = \gamma + 1$, $\beta_\gamma$ has been defined such that $\kappa_\gamma < \beta_\gamma < \kappa_{\gamma + 1} = \kappa_\eta$, and $E_{\beta_\gamma}$ is an equivalence relation. Let $H_{\gamma + 1}$ be a $\text{Coll}(\omega, \kappa_{\gamma + 1})$-generic over $M'$ filter. Since $E$ is a $\Sigma_1^1(z)$ equivalence relation, the statement ``$E$ is an equivalence relation'' is $\Pi_2^1(z)$. This statement is absolute between $M'$ and any generic extension of $M'$ by Schoenfield absoluteness. So $E$ remains a $\Sigma_1^1(z)$ equivalence relation in $M'[H_{\gamma + 1}]$. Therefore Fact \ref{Burgess theorem} holds in $M'[H_{\gamma + 1}]$ and $(\omega_1)^{M'[H_{\gamma + 1}]} = \kappa_{\gamma + 2}$. There exists some $\beta$ with $\kappa_{\gamma + 1} = \kappa_\eta < \beta < \kappa_{\gamma + 2} = \kappa_{\eta + 1}$ with $E_\beta$ an equivalence relation in $M'[H_{\gamma + 1}]$. Recall that $E_\beta$ is $\Delta_1^1(z,c)$ for any $c$ that codes $\beta$. In $M'[H_{\gamma + 1}]$, $\beta$ is countable. A real $c$ coding $E_\beta$ exists in $M'[H_{\gamma + 1}]$. $E_\beta$ is $\bDelta_1^1$ in $M'[H_{\gamma + 1}]$. Therefore, the statement ``$E_\beta$ is an equivalence relation'' is $\bPi_1^1$ in $M'[H_{\gamma + 1}]$. By Mostowski absoluteness between $M'[H_{\gamma + 1}]$ and $V$, $E_\beta$ is an equivalence relation in $V$. Let $\beta_\eta$ be this $\beta$. Finally, $E_\alpha = \bigcap_{\eta < \xi} E_{\beta_\eta}$. So $E_\alpha$ is an equivalence relation.

Now Lemma \ref{analytic equivalence stabilization of classes} can be applied to show that $[g]_E = [g]_{E_\alpha}$. This $\alpha$ is as required.
\end{proof}

\Begin{remark}{true pi11 code in generic extensions remark}
Note that the proof of the above lemma shows that if $\alpha$ is the ordinal height of the Mostowski collapse of the countable elementary substructure $M$, then $E_\alpha$ is an equivalence relation. This lemma is quite general as it only demands that $\bbP$ adds a generic real. A simpler proof can be given if $\bbP$ satisfies some additional properties:

If $\kappa$ is a cardinal, then $H_\kappa \models \mathsf{KP}$. If $1_\bbP^M \models \aleph_1 = \check \kappa$, then there is no need to collapse any $M$-cardinals to use Mostowski's absoluteness: Let $g$ be any $\bbP$-generic over $M$, Lemma \ref{Burgess theorem} gives a club subset $C$ of $\kappa = \omega_1^{M[g]}$ for which each $E_\alpha$ is an equivalence relation. $M[g]$ thinks $E_\alpha$ is a $\borel$ equivalence relation. The statement that $E_\alpha$ is a $\borel$ equivalence relation is $\bPi_1^1$ and so absolute into the real universe. As in the proof of the lemma, this can be used to show that $E_{\kappa'}$ is an equivalence relation where $\kappa'$ is the image of $\kappa$ under the Mostowski collapse of $M$. In particular, if $\bbP$ is $\aleph_1$-preserving (for example, proper), then this situation holds.

Also there are more careful versions of Lemma \ref{Burgess theorem} in which all the $E_\alpha$'s are equivalence relations which could be used to avoid this issue entirely. However, the simpler form of Lemma \ref{Burgess theorem} was used so that it could be more easily applied to the less familiar $\bPi_1^1$ setting in Lemma \ref{pi11 increasing union borel equiv}.
\end{remark}

Now returning to the setting of the main question: Suppose $E$ is a $\Sigma_1^1(z)$ equivalence relation with all classes $\bDelta_1^1$. Let $I$ be a $\sigma$-ideal with $\bbP_I$-proper. According to Lemma \ref{true pi11 code in generic extensions}, if one could find some $M \prec H_\Theta$ such that whenever $x$ is $\bbP_I$-generic over $M$, a $\Pi_1^1$ code for $[x]_E$ resides inside $M[x]$, then letting $C$ be the $I^+$ $\bDelta_1^1$ set of $\bbP_I$-generic reals over $M$, there exists some $\alpha < \omega_1$ such that $E \upharpoonright C = E_\alpha \upharpoonright C$. Hence $E \upharpoonright C$ is $\bDelta_1^1$.

A plausible candidate for the $\Pi_1^1$ code of $[x]_E$ which is an element of $M[x]$ would be some $y$ such that $M[x] \models [x]_E$ is $\Pi_1^1(y)$. However, $M[x]$ may not think $[x]_E$ is $\bDelta_1^1$. The statement ``all $E$-class are $\bDelta_1^1$'' is $\Pi_4^1(z)$. If $V$ satisfies $\bPi_4^1$-generic absoluteness, one can choose $M \prec H_\Theta$ such that some particular $\Pi_4^1(z)$ statement becomes absolute between $M$ and all its generic extensions. So in such a structure $M$, $M[x]$ will think $[x]_E$ is $\bDelta_1^1$. 

Now in $M[x]$, there is some $y$ such $M[x] \models [x]_E$ is $\Pi_1^1(y)$. In general, it is not clear if $[x]_E$ is $\Pi_1^1(y)$ in $V$. The formula ``$[x]_E$ is $\Pi_1^1(y)$'' is $\Pi_2^1(z)$. One can not use Schoenfield absoluteness between $M[x]$ (or rather it transitive collapse) and $V$ since it is not the case that $\omega_1^V \subseteq M[x]$ because $M[x]$ is countable in $V$. So what is needed is some $M \prec H_\Theta$ such that for all $\bbP_I$-generic over $M$ real $x$, a certain $\Pi_2^1(z)$ formula is absolute between $M[x]$ and $V$. The concept of universally Baireness can be used to remedy this issue.

\Begin{definition}{universally baire}
(\cite{Universally-Baire-Sets-of-Reals}) $A \subseteq \bairespace$. $A$ is universally Baire if and only if there exists $\alpha, \beta \in \text{ON}$ and trees $U$ on $\omega \times \alpha$ and $W$ on $\omega \times \beta$ such that

(1) $A = p[U]$. $\bairespace - A = p[W]$.

(2) For all $\bbP$, $1_\bbP \forces_\bbP p[\check U] \cup p[\check W] = \bairespace$. 

\noindent where $p$ of a tree denotes the projection of the tree.
\end{definition}

\Begin{fact}{tree representation sigma12 continue to represent in extension}
Suppose $A$ is a $\bSigma_2^1$ set defined by a $\bSigma_2^1$ formula $\varphi(x)$. Let $U$ and $W$ be trees witnessing that $A$ is universally Baire. Then $1_\bbP \forces_\bbP (\forall x)(\varphi(x) \Leftrightarrow x \in p[\check U])$. 
\end{fact}

\begin{proof}
See \cite{Universally-Baire-Sets-of-Reals}, page 221-222.
\end{proof}

\Begin{definition}{set of codes for pi11 sets}
Let $E$ be a $\Sigma_1^1(z)$ equivalence relation. Define the set $D$ by
$$(x,T) \in D \Leftrightarrow (T \text{ is a tree on } \omega \times \omega) \wedge (\forall y)(y \ E \ x \Leftrightarrow T^y \in \text{WF})$$
$D$ is $\Pi_2^1(z)$. 
\end{definition}

Finally, the first result showing that a positive answer follows from some strong set theoretic assumptions:

\Begin{proposition}{universally baire pi14 generic absoluteness analytic arrow borel}
Assume all $\bPi_2^1$ sets are universally Baire and $\bPi_4^1$-generic absoluteness holds. Let $I$ be a $\sigma$-ideal such that $\bbP_I$ is proper. Then $\bSigma_1^{1 \ \bDelta_1^1} \rightarrow_I \bDelta_1^1$. 
\end{proposition}

\begin{proof}
Let $E$ be a $\Sigma_1^1(z)$ equivalence relation. Since all $\bPi_2^1$ sets are universally Baire, let $U$ and $W$ be trees on $\omega \times \omega \times \alpha$ and $\omega \times \omega \times \beta$, respectively, where $\alpha, \beta \in \text{ON}$, giving the universally Baire representations for the $\Pi_2^1(z)$ set $D$ from Definition \ref{set of codes for pi11 sets}.

Suppose $B \in \bbP_I$. Using the reflection theorem, choose $\Theta$ large enough so that $B$, $\bbP_I$, $z$, $U$, and $W$ are contained in $H_\Theta$ and $H_\Theta$ satisfies $\bPi_4^1$-generic absoluteness for the statement ``$(\forall x)(\exists T)((x,T) \in D)$''. Let $M \prec H_\Theta$ be a countable elementary substructure containing $B$, $\bbP_I$, $z$, $U$, and $W$. 

By Fact \ref{properness equivalence}, let $C$ be the $I^+$ $\bDelta_1^1$ subset of $\bbP_I$-generic over $M$ reals in $B$. Let $g \in C$. Since $E$ has all classes $\bDelta_1^1$, $M$ satisfies $(\forall x)(\exists T)((x,T) \in D)$. Because $M$ has generic absoluteness for this formula, $M[g] \models (\forall x)(\exists T)((x,T) \in D)$. There exists some $T \in M[g]$ such that $M[g] \models (g,T) \in D$. 

By Fact \ref{tree representation sigma12 continue to represent in extension}, $M[g] \models (g,T) \in p[U]$. There exists $\Phi : \omega \rightarrow \alpha$ with $\Phi \in M[g]$ such that $(g,T,\Phi) \in [U]$. For each $n \in \omega$, $(g \upharpoonright n, T \upharpoonright n, \Phi \upharpoonright n) \in M$. For each $n$, $M[g] \models (g \upharpoonright n, T \upharpoonright n, \Phi \upharpoonright n) \in U$. By absoluteness, for each $n$, $M \models (g \upharpoonright n, T \upharpoonright n, \Phi \upharpoonright n) \in U$. Since $M \prec H_\Theta$, for all $n$, $(g \upharpoonright n, T \upharpoonright n, \Phi \upharpoonright n) \in U$ in the true universe $V$. Therefore, in $V$, $(g,T, \Phi) \in [U]$. $(g,T) \in p[U]$. $(g,T) \in D$. Note that $(g,T) \in D$ implies that $[g]_E$ is $\Pi_1^1(T)$. 

It has been shown that for the chosen $M$, whenever $g$ is $\bbP_I$-generic over $M$, there exists some $z \in M[g]$ such that $[g]_E$ is $\Pi_1^1(z)$. By Lemma \ref{true pi11 code in generic extensions}, there is a countable ordinal $\alpha$ such that for all $\bbP_I$-generic over $M$ reals $g$, $[g]_E = [g]_{E_\alpha}$. Hence $E \upharpoonright C = E_{\alpha} \upharpoonright C$. $E \upharpoonright C$ is $\bDelta_1^1$. 
\end{proof}

\Begin{remark}{strength of universally baire and generic absoluteness}
By \cite{Generic-Absoluteness}, Theorem 8, $\bPi_4^1$ generic absoluteness is equiconsistent with every set having a sharp and the existence of a reflecting cardinal. The proof of \cite{Generic-Absoluteness}, Theorem 8, shows that any structure satisfying $\Sigma_4^1$ generic absoluteness is closed under sharps. By \cite{Universally-Baire-Sets-of-Reals}, Theorem 3.4, all $\bPi_2^1$ sets are universally Baire is equivalent to the existence of sharps for all sets. Hence, the hypothesis of Proposition \ref{universally baire pi14 generic absoluteness analytic arrow borel} is equiconsistent with all sets having sharps and the existence of a reflecting cardinal.

Observe that $\bPi_4^1$ generic absoluteness can be avoided for those $\bSigma_1^1$ equivalence relation such that the statement ``all $E$-classes are $\bDelta_1^1$'' hold in any model of $\mathsf{ZFC}$ containing the defining parameter for $E$. 
\end{remark}

\Begin{proposition}{sharps reflecting consistency positive answer}
The consistency of $\mathsf{ZFC}$, sharps of all sets exists, and there exists a reflecting cardinal implies the consistency of $\analytic^\borel \rightarrow_I \borel$ for all $\sigma$-ideal $I$ on a Polish space such that $\bbP_I$ is proper.
\end{proposition}

\begin{proof}
See Remark \ref{strength of universally baire and generic absoluteness}.
\end{proof}

Next, a positive answer to the main question will be obtained from assumptions with weaker consistency strength. The result above illustrates all the main ideas but used stronger than necessary assumptions: $\bPi_4^1$ generic absolutness and all $\bPi_2^1$ sets are universally Baire. $\bPi_4^1$ generic absoluteness was used to preserve the statement ``all $E$-classes are $\borel$''. Below, it will be shown how sharps can be used to give a $\bPi_3^1$ statement which is equivalent. Sharps will also be used to make the statement ``all $E$-classes are $\borel$'' true in the desired generic extensions, which is more subtle than just applying Martin-Solovay absoluteness. As observed above, sharps play an important role in $\bPi_2^1$ sets being universally Baire. In the following, a much more careful analysis will be given to determine exactly which sharps are needed.

For the more optimal proof, iterable structures will be the main tools. Familar examples of iterable structures are $V$ itself when $V$ has a measurable cardinal, certain elementary substructures of $V_\Theta$ when $V$ contains a measurable cardinal, and mice that come from the existence of sharps. In the first two, the measure exists in the structure, but in the latter, the measure is external. Some references for this material are \cite{Set-Theory-Exploring}, \cite{Inner-Models-with-Large-Cardinal-Features}, and any text in inner model theory.

Let $X$ be some set. Recall a simple formulation of the statement ``$X^\sharp$ exists'' is that there is an elementary embedding $j : L[X] \rightarrow L[X]$. Another classical formulation is that there is a closed unbounded class of indiscernible (called the Silver's indiscernible) for $L[X]$. When $x \in \cantorspace$, the object $x^\sharp$ can be considered a real coding statements about indiscernibles (in a language with countably many new constant symbols to be interpreted as a countably infinite subset of indiscernibles) true in $L[x]$. Another very useful characterization of $X^\sharp$ is given by mice:

\Begin{definition}{sharps and mouse}
(See \cite{Set-Theory-Exploring}, Definition 10.18, Definition 10.30, and Definition 10.37.) Let $\SCRL = \{\dot \in, \dot E, \dot U\}$ where $\dot \in$ is a binary relation symbol, $\dot E$ is a unary predicate, and $\dot U$ is also a unary predicate. Let $X$ be a set. A $X$-mouse is a $\SCRL$-structure, $\CM = \langle J_\alpha[X], \in, X, U \rangle$ where $J_\alpha[X]$ is the $\alpha^\text{th}$ level of Jensen's fine structural hierarchy of $L[X]$, $\dot E^\CM = X$, and $\dot U^\CM = U$  with the following additional properties:

(a) $\CM$ is an amenable structure, i.e., for all $z \in J_\alpha[X]$, $z \cap X \in J_\alpha[X]$ and $z \cap U \in J_\alpha[X]$. 

(b) In the language $\{\dot \in\}$, $(J_\alpha[X], \in) \models$ ``$\mathsf{ZFC - P}$ and there is a largest cardinal''.

(c) If $\kappa$ is the largest cardinal of $(J_\alpha[X], \in)$, then $\CM \models$ $U$ is a $\kappa$-complete normal non-trivial ultrafilter on $\kappa$. 

(d) $\CM$ is iterable, i.e., every structure appearing in any putative linear iteration of $\CM$ (by $U$) is well-founded.
\\*
\\*\indent The statement $X^\sharp$ exists is also equivalent to the existence of an $X$-mouse. $X^\sharp$ will sometimes also denote the smallest $X$-mouse $\CM$ in the sense that if $\CN$ is an $X$-mouse, then there is an $\alpha$ such the $\alpha^\text{th}$ iteration $\CM_\alpha$ is $\CN$. 
\end{definition}

Under the condition that sharps of all reals exists, the statement ``all $E$-classes are $\borel$'' will be shown to be $\bPi_3^1$. This is a significant improvement since $\bPi_3^1$ generic absolutness is much easier to obtain. 

\Begin{proposition}{sharps of reals borelness of all classes pi13}
Let $E$ be a $\Sigma_1^1(z)$ equivalence relation. There is a $\Pi_3^1(z)$ formula $\varpi(v)$ in free variable $v$ such that:

Let $x \in \cantorspace$. If $(x \oplus z)^\sharp$ exists, then the statement ``$[x]_E$ is $\borel$'' is equivalent to $\varpi(x)$.

Assume for all $r \in \bairespace$, $r^\sharp$ exists. The statement ``all E-classes are $\borel$'' is equivalent to $(\forall x)\varpi(x)$. In particular, this statement is $\Pi_3^1(z)$. 
\end{proposition}

\begin{proof}
For simplicity, assume $E$ is a $\lanalytic$ equivalence relation on $\bairespace$. Let $T$ be a recursive tree on $\omega \times \omega \times \omega$ such that 
$$(x,y) \in E \Leftrightarrow T^{x,y} \text{ is illfounded}$$

Claim : Assume $x^\sharp$ exists, then
$$``V \models [x]_E \text{ is } \borel" \Leftrightarrow ``1_{\text{Coll}(\omega, < c_1)} \forces_{\text{Coll}(\omega, <c_1)} (\exists c)(c \in \text{WO} \wedge (\forall y)(\neg(y \ E \ x) \Rightarrow \text{rk}(T^{x,y}) < \text{ot}(c)))" \in x^\sharp$$
Here $c_1$ comes from $\{c_n : n \in \omega\}$, which is a collection of constant symbols used to denote indiscernibles. 

Proof of Claim: Assume $[x]_E$ is $\borel$. Then 
$$(\exists \xi < \omega_1)(\forall y)(\neg(y \ E \ x) \Rightarrow \text{rk}(T^{x,y}) < \xi)$$
Since $x^\sharp$ exists, $\omega_1$ is inaccessible in $L[x]$ and $|\mathscr{P}^{L[x]}(\text{Coll}(\omega, \xi))| = \aleph_0$. In $V$, there is a $g \subseteq \text{Coll}(\omega, \xi)$ which is $\text{Coll}(\omega, \xi)$ generic over $L[x]$. Since $g \in V$, there is a $c \in L[x][g] \subseteq V$ such that $c \in \text{WO}$ and $\text{ot}(c) = \xi$.
$$V \models (\forall y)(\neg(y \ E \ x) \Rightarrow \text{rk}(T^{x,y}) < \text{ot}(c))$$
Since this statement above is $\Pi_2^1$, Schoenfield absoluteness gives
$$L[x][g] \models (\forall y)(\neg(y \ E \ x) \Rightarrow \text{rk}(T^{x,y}) < \text{ot}(c))$$
Using the weak homogeneity of $\text{Coll}(\omega, \xi)$, 
$$L[x] \models 1_{\text{Coll}(\omega, \xi)} \forces_{\text{Coll}(\omega, \xi)} (\exists c)(c \in \text{WO} \wedge (\forall y)(\neg(y \ E \ x) \Rightarrow \text{rk}(T^{x,y}) < \text{ot}(c)))$$
The statement forced above is $\Sigma_3^1$. By upward absoluteness of $\bSigma_3^1$ statements
$$L[x] \models 1_{\text{Coll}(\omega, <\omega_1)} \forces_{\text{Coll}(\omega, <\omega_1)} (\exists c)(c \in \text{WO} \wedge (\forall y)(\neg(y \ E \ x) \Rightarrow \text{rk}(T^{x,y}) < \text{ot}(c)))$$
$$``1_{\text{Coll}(\omega, <c_1)} \forces_{\text{Coll}(\omega, <c_1)} (\exists c)(c \in \text{WO} \wedge (\forall y)(\neg(y \ E \ x) \Rightarrow \text{rk}(T^{x,y}) < \text{ot}(c)))" \in x^\sharp$$

$(\Leftarrow)$ Assume
$$``1_{\text{Coll}(\omega, <c_1)} \forces_{\text{Coll}(\omega, <c_1)} (\exists c)(c \in \text{WO} \wedge (\forall y)(\neg(y \ E \ x) \Rightarrow \text{rk}(T^{x,y}) < \text{ot}(c)))" \in x^\sharp$$
Let $\xi < \omega_1$ be a Silver indiscernible for $L[x]$. Then
$$L[x] \models 1_{\text{Coll}(\omega, <\xi)} \forces_{\text{Coll}(\omega, <\xi)} (\exists c)(c \in \text{WO} \wedge (\forall y)(\neg(y \ E \ x) \Rightarrow \text{rk}(T^{x,y}) < \text{ot}(c)))$$
Again since $\xi < \omega_1$ and $\omega_1$ is inaccessible in $L[x]$, $\mathscr{P}^{L[x]}(\text{Coll}(\omega, < \xi))$ is countable in $V$. In $V$, there exists $g \subseteq \text{Coll}(\omega, < \xi)$ which is $\text{Coll}(\omega, < \xi)$-generic over $V$. 
$$L[x][g] \models (\exists c)(c \in \text{WO} \wedge (\forall y)(\neg(y \ E \ x) \Rightarrow \text{rk}(T^{x,y}) < \text{ot}(c)))$$
Since $g \in V$, $L[x][g] \subseteq V$ and there exists a $c \in L[x][g]$ such that
$$L[x][g] \models (\forall y)(\neg(y \ E \ x) \Rightarrow \text{rk}(T^{x,y}) < \text{ot}(c))$$
This statement is $\bPi_2^1$. Since $L[x][g] \subseteq V$, Schoenfield absoluteness can be applied to give
$$V \models (\forall y)(\neg(y \ E \ x) \Rightarrow\text{rk}(T^{x,y}) < \text{ot}(c))$$
Therefore, 
$$V \models [x]_E \text{ is } \borel$$
This concludes the proof of the claim.

The statement in variables $v$ and $w$ expressing ``$w = v^\sharp$'' is $\Pi_2^1$. Therefore
$$[x]_E \text{ is } \borel$$
if and only if
$$(\forall y)((y = x^\sharp) \Rightarrow ``1_{\text{Coll}(\omega, < c_1)} \forces_{\text{Coll}(\omega, < c_1)} (\exists c)(c \in \text{WO} \wedge (\forall y)(\neg(y \ E x) \Rightarrow \text{rk}(T^{x,y}) < \text{ot}(c))) \in y)$$
The latter is $\Pi_3^1(x)$. 

Similarly
$$(\forall x)([x]_E \text{ is } \borel)$$
if and only if 
$$(\forall x)(\forall y)((y = x^\sharp) \Rightarrow ``1_{\text{Coll}(\omega, <c_1)} \forces_{\text{Coll}(\omega, <c_1)} (\exists c)(c \in \text{WO} \wedge (\forall y)(\neg(y \ E \ x) \Rightarrow \text{rk}(T^{x,y}) < \text{ot}(c)))" \in y)$$
The latter is $\Pi_3^1$. 

Let $\varpi(v)$ be the statement:
$$(\forall y)(y = v^\sharp \Rightarrow ``1_{\text{Coll}(\omega, < c_1)} \forces_{\text{Coll}(\omega, <c_1)} (\exists c)(c \in \text{WO} \wedge (\forall y)(\neg(y \ E \ v) \Rightarrow \text{rk}(T^{v,y}) < \text{ot}(c)))" \in y)$$
By the above results, this works.
\end{proof}

So assuming for all $x \in \cantorspace$, $x^\sharp$ exists, the statement ``all $E$-classes are $\borel$'' is $\bPi_3^1$. Below some conditions for $\bPi_3^1$ generic absoluteness will be explored. However, there is still a subtle point to be noted. Assuming all sharps of reals exists and generic $\bPi_3^1$ absoluteness holds for a forcing $\bbP$. Let $G \subseteq \bbP$ be $\bbP$-generic over $V$. Then in $V[G]$, the statement $(\forall x)(\varpi(x))$ remains true by $\bPi_3^1$ absoluteness. But if $V[G]$ does not satisfy all sharps of reals exists, then it may not be true that $(\forall v)(\varpi(v))$ is equivalent to the statement ``all $E$-classes are $\borel$''. For the main question in the case of  $\bbP_I$, only the fact that $[\dot x_\text{gen}]_E$ is $\borel$ will be of any concern. In $V[G]$, one has $(\forall x)(\varpi(x))$. In particular, if $g$ is the generic real added by $G$, then $\varpi(g)$ holds. If $(g \oplus z)^\sharp$ exists, then Proposition \ref{sharps of reals borelness of all classes pi13} implies ``$[g]_E$ is $\borel$'' is equivalent to $\varpi(g)$. Hence $[g]_E$ is $\borel$ in $V[G]$. The following is a situation (applicable later) for which $(g \oplus z)^\sharp$ exists. 

\Begin{fact}{forcing extension has sharps}
Let $A$ be a set. Suppose $A^\sharp$ exists and $j : L[A] \rightarrow L[A]$ is a nontrivial elementary embedding. Suppose $\bbP \in L[A]$ is a forcing in $(V_{\text{crit}(j)})^{L[A]}$. Suppose $G \subseteq \bbP$ is generic over $V$ (or just $L[A]$), then $\langle A, G\rangle^\sharp$ exists in $V[G]$. 
\end{fact}

\begin{proof}
Since $\bbP \in (V_{\text{crit}(j)})^{L[A]}$, $j'' \ \text{tc}(\{\bbP\}) = \text{tc}(\{\bbP\})$. Define the lift $\tilde j : L[\langle A, G \rangle] \rightarrow L[\langle A, G \rangle]$ by
$$\tilde j (\tau[G]) = j(\tau)[G]$$ 
By the usual arguments, $\tilde j$ is a nontrivial elementary embedding definable in $V[G]$. Hence $\langle A, G \rangle^\sharp$ exists.
\end{proof}

Next, a few more basic properties of iterable structures:

\Begin{fact}{structure elementary embed onto iterable structure is iterable}
Let $\SCRL = \{\dot \in, \dot U\}$. Suppose $\CN = (N, \in, V)$ is an iterable structure. Suppose $\CM = (M, \in, U)$ is a $\SCRL$-structure such that there exists an $\SCRL$-elementary embedding $j: \CM \rightarrow \CN$. Then $\CM$ is iterable. 
\end{fact}

\begin{proof}
See \cite{Inner-Models-with-Large-Cardinal-Features}, Lemma 18 for a proof.
\end{proof}

\Begin{fact}{version of iterable in elementary substructure is also iterable}
Let $M \prec H_\Theta$ be a countable elementary substructure where $\Theta$ is some sufficienty large cardinal. Let $\mathcal{U}$ be an iterable structure and $\mathcal{U} \in M$. Let $\mathcal{U}^M = \mathcal{U} \cap M$. Then $\mathcal{U}^M$ is iterable. 
\end{fact}

\begin{proof}
Let $\varphi$ be some $\SCRL = \{\dot \in, \dot U\}$ sentence. For any $x \in \mathcal{U} \cap M$, $\mathcal{U}^M \models \varphi(x)$ if and only if $M \models \mathcal{U} \models \varphi(x)$. Since $M \prec H_\Theta$, if and only if $V \models \mathcal{U} \models \varphi(x)$. Hence $\mathcal{U}^M \prec \mathcal{U}$ as an $\SCRL$-structure. By Fact \ref{structure elementary embed onto iterable structure is iterable}, $\mathcal{U}^M$ is iterable.
\end{proof}

As mentioned above, it is not possible in general to claim that $\bPi_2^1$ statements are absolute between a countable model $M$ and the universe $V$ since Schoenfield absoluteness can not be applied when it is not the case that $\omega_1^V \subseteq M$. However, $\omega_1$-iterable structures can be used to solve this problem by applying Schoenfield absoluteness in the $\omega_1$ iteration.

\Begin{fact}{mouse pi12 correctness}
Let $X$ be a set. Suppose $\CM = (J_\alpha[X], \in, X, U)$ is an $X$-mouse. Then $J_{\alpha}[X]$ is $\Pi_2^1$-correct, that is, if $\varphi$ is a $\bPi_2^1$ sentence with parameters in $J_\alpha[X]$, then $J_{\alpha}[X] \models \varphi$ if and only if $V \models \varphi$. 

Let $\kappa$ be the largest cardinal of $(J_\alpha[X], \in)$. Suppose $\bbP \in J_\alpha[X]$ is a forcing such that $(J_\alpha[X], \in) \models \bbP \in V_\kappa$. Then $J_{\alpha}[X]$ is $\bbP$-generically $\Pi_2^1$-correct, that is, for all $G \subseteq \bbP$ which is $\bbP$-generic over $J_\alpha[X]$ and $G \in V$, and any $\bPi_2^1$-formula coded in $J_\alpha[X][G]$, $J_\alpha[X][G] \models \varphi$ if and only if $V \models \varphi$. 
\end{fact}

\begin{proof}
Let $\CM_0 = \CM$. Let $j_{0, \omega_1} : \CM_0 \rightarrow \CM_{\omega_1}$ denote the $\omega_1$-iteration. $\CM_{\omega_1}$ is well-founded, so let $\CM_{\omega_1} = (J_\beta[X], \in, X, U_{\omega_1})$. Note that $\beta \geq \omega_1$. By \cite{Set-Theory-Exploring}, Lemma 10.21 (d), $j_{0,\omega_1}$ is a full ($\Sigma_\omega$) elementary embedding in the language $\{\dot\in, \dot E\}$. So if $\varphi$ is a $\bPi_2^1$ statement with parameter in $J_\alpha[X]$ then $J_\alpha[X] \models \varphi$ if and only if $J_\beta[X] \models \varphi$. Since $\omega_1^V \subseteq J_\beta[X]$ , $J_\beta[X] \models \varphi$ if and only if $V \models \varphi$. 

For the second statement: Since $(J_\alpha[X], \in) \models \bbP \in V_\kappa$, $j_{0,\omega_1}$ does not move any elements in the transitive closure of $\bbP$. Also no new subsets of $\bbP$ appears in $J_\beta[X]$. Thus if $G$ is $\bbP$-generic over $J_\alpha[X]$, then $G$ is $\bbP$-generic over $J_\beta[X]$. Lift the elementary embedding  $j_{0,\omega_1} : J_\alpha[X] \rightarrow J_\beta[X]$ to $\tilde j_{0,\omega_1} : J_\alpha[X][G] \rightarrow J_\beta[X][G]$ in the usual way: if $\tau \in J_\alpha[X]^\bbP$, then 
$$\tilde j_{0,\omega_1}(\tau[G]) = j_{0,\omega_1}(\tau)[G]$$  
$\tilde j_{0,\omega_1}$ is a well-defined elementary embedding. Let $\varphi$ be a $\Pi_2^1$ formula coded in $J_\alpha[X][G]$. Using this elementary embedding, $J_\alpha[X][G] \models \varphi$ if and only if $J_\beta[X][G] \models \varphi$. Since $\omega_1^V \subseteq J_\beta[X][G]$ and using Schoenfield absoluteness, $J_\beta[X][G] \models \varphi$ if and only if $V \models \varphi$. 
\end{proof}

\Begin{fact}{name associated statement}
Let $\bbP$ be a forcing. Suppose $\varphi(v)$ is a formula with one free variable. By fullness or the maximality principle (see \cite{Set-Theory-Kunen} Theorem IV.7.1), there exists a name $\tau^\bbP_\varphi$ such that $1_\bbP \forces_{\bbP} (\exists v)(\varphi(v))$ if and only if $1_\bbP \forces_{\bbP} \varphi(\tau_\varphi^\bbP)$. If $\tau_\varphi^\bbP$ is a name for a real, one may assume that it is a nice name for a real. 
\end{fact}

\Begin{fact}{pi31 absoluteness for particular forcing}
Consider the $\bSigma_3^1$ sentence $(\exists v)(\varphi(v))$ where $\varphi(v)$ is $\Pi_2^1$. If $\langle \bbP, \tau_\varphi^\bbP\rangle^\sharp$ exists, then $(\exists v)(\varphi(v))$ is absolute between the ground model and $\bbP$-extensions.
\end{fact}

\begin{proof}
This is originally proved using the Martin-Solovay tree, which were implicit in \cite{A-Basis-Theorem-for-Sigma13-Sets-of-Reals}. The proof from \cite{Projective-Well-Orderings-of-the-Reals} Theorem 3 is sketched below to make explicit what sharps are necessary.

Suppose $1_\bbP \forces_\bbP^V (\exists v)\varphi(v)$. Then $1_\bbP \forces_{\bbP}^V \varphi(\tau_\varphi^\bbP)$. 

Note that $\bbP \in \langle \bbP, \tau_\varphi^\bbP\rangle^\sharp$ and $\langle \bbP, \tau_\varphi^\bbP\rangle^\sharp \models \bbP \in V_\kappa$, where $\kappa$ is the largest cardinal of $\langle \bbP, \tau_\varphi^\bbP\rangle^\sharp$.  

Using some standard way of coding, let $T$ be a tree of attempts to build a tuple $(\CM, \bbQ, H, y, j)$ with the following properties:

(1) $\CM$ is a countable structure satisfying (a), (b), and (c) from Definition \ref{sharps and mouse}.

(2) $\bbQ \in \CM$ is a forcing.

(3) $H$ is $\bbQ$-generic over $\CM$. 

(4) $y$ is a real in $\CM[H]$ and $\CM[H] \models \varphi(y)$. 

(5) $j : \CM \rightarrow \langle \bbP, \tau_\varphi^\bbP\rangle^\sharp$ is an elementary embedding in the language $\{\dot\in, \dot E, \dot U\}$ with $j(\bbQ) = \bbP$. 

\noindent Let $G$ be an arbitrary $\bbP$-generic over $V$. Since $V[G] \models (\exists v) \varphi(v)$, $V[G] \models \varphi(\tau_\varphi^\bbP[G])$. By the downward absoluteness of $\Pi_2^1$ statements (which follows from Mostowski absoluteness), $\langle \bbP, \tau_\varphi^\bbP\rangle^\sharp[G] \models \varphi(\tau_\varphi^\bbP[G])$. By Downward-Lowenheim-Skolem, let $\CN$ be a countable $\{\dot \in, \dot E, \dot U\}$ elementary substructure of $\langle \bbP, \tau_\varphi^\bbP\rangle^\sharp$ containing $\bbP$ and $\tau_\varphi^\bbP$. Let $\CM$ be the Mostowski collapse of $\CN$, and $j : \CM \rightarrow \langle \bbP, \tau_\varphi^\bbP\rangle^\sharp$ be the induced elementary embedding. Let $\bbQ = j^{-1}(\bbP)$. Let $H = j^{-1}[G]$. Let $y = j^{-1}(\tau_\varphi^\bbP)[H]$. So in $V[G]$, $(\CM, \bbQ, H, y, j)$ is a path through $T$. 

Therefore, in $V[G]$, the tree $T$ is illfounded. Hence it is illfounded in $V$ by $\Delta_1$-absoluteness. In $V$, let $(\CM, \bbQ, H, y, j)$ be such a path. By Fact \ref{structure elementary embed onto iterable structure is iterable}, $\CM$ is iterable. By Fact \ref{mouse pi12 correctness}, $\CM[H] \models \varphi(y)$ implies $V \models \varphi(y)$. This establishes that $(\exists v)\varphi(v)$ is downward absolute from $V[G]$ to $V$. This completes the proof.
\end{proof}

\Begin{definition}{name associated with all classes borel}
Let $I$ be a $\sigma$-ideal on a Polish space $X$ such that $\bbP_I$ is proper. If $\varpi$ is the formula from Proposition \ref{sharps of reals borelness of all classes pi13}, then $\neg(\forall v)\varpi(v)$ is $\bSigma_3^1$ and can be written as $(\exists v)\zeta(v)$ where $\zeta$ is $\bPi_2^1$. Let $\mu_E^I$ be $\tau_{\zeta}^{\bbP_I}$ from Fact \ref{name associated statement}.
\end{definition}

\Begin{definition}{name associated all classes pi11}
Let $I$ be a $\sigma$-ideal on $\bairespace$ such that $\bbP_I$ is proper.  Consider the formula ``$(\exists y)([\dot x_\text{gen}]_E \text{ is } \Pi_1^1(y))$''. Write it as $(\exists y)\psi(y)$. By Fact \ref{name associated statement}, let $\sigma_E^{I}$ be $\tau_\psi^{\bbP_I}$. 
\end{definition}

\Begin{definition}{associated triple for consistency main question}
Suppose $I$ is a $\sigma$-ideal on $\bairespace$ such that $\bbP_I$ is proper. Let $E \in \analytic^\borel$. Define $\chi_E^I = \langle \bbP_I, \mu_E^{I}, \sigma_E^{I}\rangle$. 
\end{definition}

The following result gives a positive answer to the main question for $\analytic^\borel$ using sharps for some small sets. 

\Begin{theorem}{sharps of reals and positive answer main question}
Suppose $I$ is a $\sigma$ ideal on $\bairespace$ such that $\bbP_I$ is proper. If for all $x \in \cantorspace$, $x^\sharp$ exists and $(\chi_E^I)^\sharp$ exists for all $E \in \analytic^\borel$, then $\analytic^\borel \rightarrow_I \borel$. 
\end{theorem}

\begin{proof}
Let $\Theta$ be sufficiently large and $M \prec H_\Theta$ is countable elementary with $(\chi_E^I)^\sharp \in M$. Note that $M \models (\chi_E^I)^\sharp$ exists and for all $x \in \cantorspace$, $x^\sharp$ exists. 

First to show that $(\forall v)(\varpi(v))$ is a $\bbP_I$-generically absolute for $M$: Since $M$ satisfies all sharps of reals exists, Proposition \ref{sharps of reals borelness of all classes pi13} implies ``all $E$-classes are $\borel$'' is equivalent to $(\forall v)(\varpi(v))$. The latter is $\bPi_3^1$ and so its negation is $\bSigma_3^1$. Since $M \models \langle \bbP_I, \mu_E^I\rangle^\sharp$ exists, Fact \ref{pi31 absoluteness for particular forcing} implies the statement, $(\forall v)(\varpi(v))$ is absolute between $M$ and $\bbP_I$ extensions of $M$. Since $M$ satisfies ``all $E$-classes are $\borel$'' and $M$ satisfies all sharps of reals exists, $M$ satisfies $(\forall v)(\varpi(v))$. Therefore, all $\bbP_I$ extensions of $M$ satisfy the formula $(\forall v)(\varpi(v))$. 

Since $\bbP_I^\sharp$ exists, there exists a $j : L[\bbP_I] \rightarrow L[\bbP_I]$ with $\bbP_I \in (V_{\text{crit}(j)})^{L[\bbP]}$. Therefore, Fact \ref{forcing extension has sharps} implies $1_{\bbP_I} \forces_{\bbP_I} \dot x_\text{gen}^\sharp \text{ exists }$. So $M \models 1_{\bbP_I} \forces_{\bbP_I} \dot x_\text{gen}^\sharp \text{ exists }$. By Proposition \ref{sharps of reals borelness of all classes pi13}, 
$$M \models 1_{\bbP_I} \forces_{\bbP_I} ([\dot x_\text{gen}]_E \text{ is } \borel) \Leftrightarrow \varpi(\dot x_\text{gen})$$
Since all $\bbP_I$ extensions of $M$ satisfy $(\forall v)(\varpi(v))$, all $\bbP_I$ extensions of $M$ satisfy $[\dot x_\text{gen}]_E$ is $\borel$. 

By the result of the previous paragraph, $1_{\bbP_I} \forces_{\bbP_I}^M (\exists y)\psi(y)$, where $\psi$ is the formula from Definition \ref{name associated all classes pi11}. Therefore, $1_{\bbP_I} \forces_{\bbP_I}^M \psi(\sigma_E^I)$. Note that $\psi(y)$ is actually $\psi'(\dot x_\text{gen}, y)$, where $\psi'$ is $\bPi_2^1$ with parameters from $M$ asserting that $[\dot x_\text{gen}]$ is $\Pi_1^1(y)$. Since $\dot x_\text{gen}$ is constructible from $\bbP_I$, the existence of $\langle \bbP_I, \dot x_\text{gen}, \sigma_E^I \rangle^\sharp$ follows from the existence of $\langle \bbP_I, \sigma_E^I\rangle^\sharp$. Applying the downward absoluteness of $\bPi_2^1$ statement from $M[x]$ to $(\langle \bbP_I, \sigma_E^I\rangle^\sharp)^M[x]$ (where $x$ is any $\bbP_I$-generic over $M$ real) gives $(\langle \bbP_I, \sigma_E^I\rangle^\sharp)^M[x] \models \psi'(x,\sigma_E^I[x])$. By Fact \ref {version of iterable in elementary substructure is also iterable}, $(\langle \bbP_I, \sigma_E^I\rangle^\sharp)^M$ is still iterable. Applying Fact \ref{mouse pi12 correctness} (generic $\Pi_2^1$-correctness) to $(\langle \bbP_I, \sigma_E^I\rangle^\sharp)^M[x]$ and $V$, one has $V \models \psi'(x, \sigma_E^I[x])$ where $x$ is any $\bbP_I$-generic over $M$. So it has been shown that $M[x] \models [x]_E$ is $\Pi_1^1(\sigma_E^I[x])$ and $V \models [x]_E$ is $\Pi_1^1(\sigma_E^I[x])$. 

Lemma \ref{true pi11 code in generic extensions} implies that there is some countable $\alpha$ such that for all $x$ which is $\bbP_I$-generic over $M$, $[x]_{E_\alpha} = [x]_E$. Therefore if $B$ is an arbitrary $I^+$ $\borel$ subset and $C$ is the $I^+$ $\borel$ set of $\bbP_I$-generic over $M$ reals in $B$, then $E \upharpoonright C = E_\alpha \upharpoonright C$. $\{E\} \rightarrow_I \borel$. 
\end{proof}

Since $\bbP_I$ is a collection of subsets of $\bairespace$ and $\mu_E^I$ and $\sigma_E^I$ can be taken to be nice names for reals, $\chi_E^I$ is an element of $H_{(2^{\aleph_0})^+}$. Therefore, if there is a measurable or a Ramsey cardinal, then $(\chi_E^I)^\sharp$ will exist. 

\Begin{corollary}{sharp hereditarily continuum plus canonicalization analytic}
If $z^\sharp$ exists for all $z \in H_{(2^{\aleph_0})^+}$, then $\analytic^\borel \rightarrow_I \borel$ for all $\sigma$-ideal $I$ such that $\bbP_I$ is proper.
\end{corollary}

\Begin{corollary}{ramsey measurable imply positive answer main question}
If there exists a Ramsey cardinal, then $\analytic^\borel \rightarrow_I \borel$ for all $\sigma$-ideal $I$ such that $\bbP_I$ is proper.
\end{corollary}

\section{$\bSigma_1^1$ Equivalence Relations and Some Ideals}\label{bsigma Equivalence Relations and Some Ideals}

Some partial results about the main question for $\analytic$ equivalence relations with all classes $\borel$ will be provided in this and the next section. These are proved using various different techniques and different set theoretic assumptions (usually of lower consistency strength than the full answer of the previous section). These results may be useful in understanding what combination of universes, $\analytic$ equivalence relations, and $\sigma$-ideals can not be used to demonstate the consistency of a negative answer to the main question.

In this section, the focus will be on the main question in the case of some classical ideals $I$ with $\bbP_I$ proper.

\Begin{definition}{countable ideal}
Let $X$ be a Polish space. Let $I_\text{countable} := \{A \subseteq X : |A| \leq \aleph_0\}$. $\bbP_{I_\text{countable}}$ is forcing equivalent to Sacks forcing.
\end{definition}

\Begin{proposition}{analytic equiv countable ideal arrow borel}
$\bSigma_1^1 \rightarrow_{I_\text{countable}} \bDelta_1^1$. 
\end{proposition}

\begin{proof}
Let $E$ be any $\bSigma_1^1$ equivalence relation. (Note there is no condition on the classes being $\bDelta_1^1$ for this proposition.) Let $B$ be $I_\text{countable}^+$ $\bDelta_1^1$ set, i.e. an uncountable $\bDelta_1^1$ set. 

Suppose there is some $x \in B$ such that $[x]_E \cap B$ is uncountable. The perfect set property for the $\bSigma_1^1$ set $[x]_E \cap B$ implies that $[x]_E \cap B$ has a perfect subset $C$. Then $E \upharpoonright C = \text{ev} \upharpoonright C$. So $\{E\} \rightarrow_{I_\text{countable}} \bDelta_1^1$. 

Otherwise, $[x]_E \cap B$ is countable for all $x \in B$. $E \upharpoonright B$ is a $\bSigma_1^1$ equivalence relation with all classes countable. Then $\{E\} \rightarrow_I \bDelta_1^1$ follows from Fact \ref{all classes countable arrow borel}.
\end{proof}

\Begin{definition}{E0 silver ideal}
Let $I_{E_0}$ denote the $\sigma$-ideal $\sigma$-generated by the $\bDelta_1^1$ sets on which $E_0$ is smooth.
\end{definition}

\Begin{fact}{PIE0 equivalent silver forcing}
$\bbP_{I_{E_0}}$ is forcing equivalent to Prikry-Silver forcing. Hence $\bbP_{I_{E_0}}$ is a proper forcing.
\end{fact}

\begin{proof}
See \cite{Descriptive-Set-Theory-and-Definable-Forcing}, Lemma 2.3.37.
\end{proof}

\Begin{fact}{analytic equiv IE0 arrow id ev E0}
$\bSigma_1^1 \rightarrow_{I_{E_0}} \{\text{id}, \text{ev}, E_0\}$
\end{fact}

\begin{proof}
See \cite{Canonical-Ramsey-Theory-on-Polish-Spaces}, Theorem 7.1.1.
\end{proof}

\Begin{corollary}{analytic equiv IEO arrow borel}
$\bSigma_1^1 \rightarrow_{I_{E_0}} \bDelta_1^1$. 
\end{corollary}

\begin{proof}
This follows immediately from Fact \ref{analytic equiv IE0 arrow id ev E0} since $\text{id}$, $\text{ev}$, and $E_0$ are all $\bDelta_1^1$ equivalence relations.
\end{proof}

\Begin{definition}{meager and null ideal}
Let $I_\text{meager}$ be the $\sigma$-ideal $\sigma$-generated by the meager subsets of $\bairespace$ (or more generally any Polish space).

Let $I_\text{null}$ be the $\sigma$-ideal $\sigma$-generated by the Lebesgue null subsets of $\bairespace$. 
\end{definition}

Kechris communicated to the author the following results concerning the meager ideal. Define a set to be $I_\text{meager}$ measurable if and only if that set has the Baire property. Define a set to be $I_\text{null}$ measurable if and only if that set is Lebesgue measurable. First, a well-known result on the additivity of the meager and null ideal under certain types of unions. 

\Begin{fact}{additivity meager null prewellordering}
Let $I$ be $I_\text{meager}$ or $I_\text{null}$. Let $\{A_\eta\}_{\eta < \xi}$ be a sequence of sets in $I$. Define a prewellordering $\sqsubseteq$ on $\bigcup_{\eta < \xi} A_\eta$ by: $x \sqsubseteq y$ if and only if the least $\eta$ such that $x \in A_\eta$ is less than or equal to the least $\eta$ such that $y \in A_\eta$. If $\sqsubseteq$ is $I$-measurable (with the version of I defined on the product space), then $\bigcup_{\eta < \xi} A_\eta$ is in $I$.
\end{fact}

\begin{proof}
See \cite{Measure-and-Category-in-Effective}, Proposition 1.5.1 for a proof.
\end{proof}

\Begin{theorem}{pi13 measurable arrow borel}
Let $I$ be $I_\text{meager}$ or $I_\text{null}$. If all $\bPi_3^1$ sets are $I$ measurable, then $\bSigma_1^{1 \ \bDelta_1^1} \rightarrow_I \bDelta_1^1$. Moreover, if $E$ is a $\bSigma_1^1$ equivalence relation with all classes $\bDelta_1^1$ and $B$ is $I^+$ $\bDelta_1^1$, then there exists a $I^+$ $\bDelta_1^1$ $C \subseteq B$ with $B \ \backslash \ C \in I$ and $E \upharpoonright C$ is a $\bDelta_1^1$ equivalence relation. 
\end{theorem}

\begin{proof}
(Kechris) For simplicity, assume $E$ is an equivalence relation on the $\bairespace$. Let $B \subseteq \bairespace$ be $I^+$ $\bDelta_1^1$. For simplicity, assume $B = \bairespace$. $(\bairespace)^2 \ \backslash \ E$ is a $\bPi_1^1$ set. Let $T$ be a tree on $\omega \times \omega \times \omega$ such that 
$$\neg (x \ E \ y) \Leftrightarrow T^{(x,y)} \text{ is well-founded}$$
For each $\alpha < \omega_1$, let $A_\alpha = \{x : (\forall y)((x,y) \notin E \Rightarrow \text{rk}(T^{(x,y)}) < \alpha)\}$. First, the claim is that for all $x \in \bairespace$, there exists some $\alpha < \omega_1$ with $x \in A_\alpha$: To see this, fix $x$ and let $L = \{(x,y) : y \notin [x]_E\}$. Since $[x]_E$ is $\bDelta_1^1$, $L$ is $\bDelta_1^1$. $L \subseteq (\bairespace)^2 \ \backslash \ E$. By the boundedness theorem, there exists some $\alpha < \omega_1$ such that $\sup \{\text{rk}(T^{(x,y)}) : (x,y) \in L\} < \alpha$. $x \in A_\alpha$. It has been shown that $\bairespace = \bigcup_{\alpha < \omega_1} A_\alpha$. 

The next claim is that there exists some $\alpha < \omega_1$ such that $A_\alpha$ is $I^+$: Suppose that for all $\alpha < \omega_1$, $A_\alpha \in I$. Note that there is a $\bPi_2^1$ formula $\Phi(x,c)$ (using the tree $T$ as a parameter) such that if $c \in \text{WO}$, then
$$\Phi(x,c) \Leftrightarrow (\forall y)(\text{rk}(T^{(x,y)}) < \text{ot}(c))$$
Define $\sqsubseteq$ using the sequence $\{A_\alpha : \alpha < \omega_1\}$. Then
$$x \sqsubseteq y \Leftrightarrow (\forall c)(c \in \text{WO} \Rightarrow (\Phi(y,c) \Rightarrow \Phi(x,c)))$$
$\sqsubseteq$ is $\bPi_3^1$ on $\bairespace \times \bairespace$. By Fact \ref{additivity meager null prewellordering}, $\bairespace = \bigcup_{\alpha < \omega_1} A_\alpha$ is in $I$. Contradiction.

Choose an $\alpha < \omega_1$ such that $A_\alpha$ is $I^+$. Since $\bPi_2^1$ sets are $I$-measurable, let $C$ be $\bDelta_1^1$ $I^+$ such that $A_\alpha \triangle C \in I$. Thus $C \ \backslash \ A_{\alpha} \in I$. Since $I$ is the $\sigma$-generated by certain $\bDelta_1^1$ sets, there exists a $\bDelta_1^1$ set $D \in I$ such that $C \ \backslash \ A_\alpha \subseteq D$. Let $B_0 = C \ \backslash \ D$. Note that $B_0$ is $I^+$ $\bDelta_1^1$ and $B_0 \subseteq A_\alpha$. 

Now suppose $\xi < \omega_1$ and a sequence $\{B_\eta : \eta < \xi\}$ of $\bDelta_1^1$ $I^+$ sets has been defined with the property that if $\eta_1 \neq \eta_2$ then $B_{\eta_1} \cap B_{\eta_2} \in I$. Let $K_\xi = \bigcup_{\eta < \xi} B_\eta$. Define $A_\alpha^\xi = A_\alpha \ \backslash \ K_\xi$. $\bairespace \ \backslash \ K_\xi = \bigcup_{\alpha < \omega_1} A_\alpha^\xi$. If $\bairespace \ \backslash \ K_\xi$ is $I^+$, then repeating the above procedure produces some $I^+$ $\bDelta_1^1$ $B_\xi$ with the property that for all $\eta < \xi$, $B_\eta \cap B_\xi \in I$ and for some $\alpha < \omega_1$, $B_\xi \subseteq A_\alpha^\xi \subseteq A_\alpha$. 

Observe that for some $\xi < \omega_1$, $\bairespace \ \backslash \ K_\xi$ must be in $I$. This is because otherwise the construction succeeds in producing an antichain $\{B_\eta : \eta < \omega_1\}$ of cardinality $\aleph_1$ in $\bbP_I$. However, $\bbP_I$ has the $\aleph_1$-chain condition. Contradiction.

So choose $\xi$ such that $\bairespace \ \backslash \ K_\xi \in I$. By construction, for each $\eta < \xi$, there is some $\alpha_\eta < \omega_1$ such that $B_\eta \subseteq A_{\alpha_\eta}^\eta \subseteq A_{\alpha_\eta}$. Since $\xi < \omega_1$, there is some $\mu < \omega_1$ such that $\sup\{\alpha_\eta : \eta < \xi\} < \mu$. Then $K_\xi = \bigcup_{\eta < \xi} A_{\alpha_\eta} \subseteq A_\mu$. Hence for all $x, y \in K_\xi$,
$$x \ E \ y \Leftrightarrow \text{rk}(T^{(x,y)}) \geq \mu$$
So $K_\xi$ is $I^+$ $\bDelta_1^1$ with $\bairespace \ \backslash \ K_\xi \in I$ and $E \upharpoonright K_\xi$ is a $\bDelta_1^1$ equivalence relation.
\end{proof}

\Begin{theorem}{large cardinal consistency imply consistency meager random arrow borel} 
The consistency of $\mathsf{ZFC}$ implies the consistency of $\mathsf{ZFC}$ and $\bSigma_1^{1 \ \bDelta_1^1} \rightarrow_{I_\text{meager}} \bDelta_1^1$. 

\noindent The consistency of $\mathsf{ZFC + Inaccessible \ Cardinal}$ implies the consistency of $\mathsf{ZFC}$ and $\bSigma_1^{1 \ \bDelta_1^1} \rightarrow_{I_{\text{null}}} \bDelta_1^1$. 
\end{theorem}

\begin{proof}
By \cite{Can-You-Take-Solovays-Inaccessible-Away}, from a model $\mathsf{ZFC}$, one can obtain a model of $\mathsf{ZFC}$ in which all $\text{OD}_{\bairespace}$ subsets of $\bairespace$ have the Baire property.

By \cite{A-Model-of-Set-Theory-in-which-Every-Set}, from a model of $\mathsf{ZFC}$ with an inaccessible cardinal, one can obtain a model of $\mathsf{ZFC}$ in which all $\text{OD}_\bairespace$ subsets of $\bairespace$ are Lebesgue measurable.

Then both results follow from Theorem \ref{pi13 measurable arrow borel}.
\end{proof}

Let $\kappa$ be an inaccessible cardinal. $\text{Coll}(\omega, <\kappa)$ denotes the L\'{e}vy collapse of $\kappa$ to $\omega_1$. Since the generic extension of the L\'{e}vy collapse of an inaccessible to $\omega_1$ (and the related Solovay's model) appears often in descriptive set theory, the following is worth mentioning:

\Begin{corollary}{levy collapse meager null arrow borel}
Let $\kappa$ be an inaccessible cardinal in $V$. Let $G \subseteq \text{Coll}(\omega, < \kappa)$ be $\text{Coll}(\omega,<\kappa)$-generic over $V$. Then in $V[G]$, $\bSigma_1^{1 \ \bDelta_1^1} \rightarrow_{I_\text{meager}} \bDelta_1^1$ and $\bSigma_1^{1 \ \bDelta_1^1} \rightarrow_{I_\text{null}} \bDelta_1^1$. 
\end{corollary}

\begin{proof}
\cite{A-Model-of-Set-Theory-in-which-Every-Set} shows that in this model, all $\text{OD}_\bairespace$ subsets of $\bairespace$ have the Baire property and is Lebesgue measurable. As above, the result follows from Theorem \ref{pi13 measurable arrow borel}.
\end{proof}

\cite{Can-You-Take-Solovays-Inaccessible-Away} shows that the existence of an inaccessible cardinal and the statement that all $\bPi_3^1$ sets are Lebesgue measurable are equiconsistent.

To show that the above statement even for $I_\text{null}$ is consistent relative to $\mathsf{ZFC}$ will require a slight modification of the above proof using a different set theoretic assumption.

\Begin{definition}{covering number}
Let $I$ be a $\sigma$-ideal on a Polish space $X$. $\text{cov}(I)$ is the smallest cardinal $\kappa$ such that there exists a set $U \subseteq I$ with $\bigcup U = X$ and $|U| = \kappa$.  
\end{definition}

\Begin{proposition}{pi12 measurable cov(I) meager null arrow borel}
Let $I$ be $I_\text{meager}$ or $I_\text{null}$. If all $\bPi_2^1$ sets are $I$-measurable and $\text{cov}(I) > \aleph_1$, then $\bSigma_1^{1\ \bDelta_1^1} \rightarrow_I \bDelta_1^1$. 
\end{proposition}

\begin{proof}
The proof is similar to Theorem \ref{pi13 measurable arrow borel}. In this case, one can conclude that for some $\alpha < \omega_1$, $A_\alpha$ is $I^+$ from the fact that $\text{cov}(I) > \aleph_1$ and $\bairespace = \bigcup_{\alpha < \omega_1} A_\alpha$. The $I$-measurability of $\bPi_2^1$ sets is needed to find some $C \subseteq A_\alpha$ which is $I^+$ $\bDelta_1^1$. 
\end{proof}

\Begin{fact}{MA not CH pi12 measurable cov(I) condition}
Let $I$ be $I_\text{meager}$ or $I_\text{null}$. $\mathsf{MA + \neg CH}$ implies all $\bPi_2^1$ sets are $I$-measurable and $\text{cov}(I) = 2^{\aleph_0} > \aleph_1$. 
\end{fact}

\begin{proof}
See \cite{Internal-Cohen-Extensions}.
\end{proof}

\Begin{theorem}{consistency ZFC consistency null meager arrow borel}
The consistency of $\mathsf{ZFC}$ implies the consistency of $\mathsf{ZFC}$ and $\bSigma_1^{1\ \bDelta_1^1} \rightarrow_I \bDelta_1^1$ where $I$ is $I_\text{meager}$ or $I_\text{null}$. 
\end{theorem}

\begin{proof}
The consistency of $\mathsf{ZFC}$ implies the consistency of $\mathsf{ZFC + MA + \neg CH}$ by a well-known iterated forcing argument. 
\end{proof}

\section{Thin $\bSigma_1^1$ Equivalence Relations}\label{Thin bSigma11 Equivalence Relations}
\Begin{definition}{thin equivalence relation}
An equivalence relation $E$ on a Polish space $X$ is thin if and only if there does not exists a perfect set of pairwise $E$-inequivalent elements.

Let $\bSigma_1^{1 \ \text{thin}}$ denote the class of thin $\bSigma_1^1$ equivalence relations defined on $\bDelta_1^1$ subsets of Polish spaces.

Let $\bSigma_1^{1 \ \text{thin} \ \bDelta_1^1}$ denote the class of thin $\bSigma_1^1$ equivalence relations with all classes $\bDelta_1^1$ defined on $\bDelta_1^1$ subsets of Polish spaces.
\end{definition}

\Begin{fact}{number of classes of thin analytic equivalence relation}
Suppose $E$ is a thin $\bSigma_1^1$ equivalence relation, then $E$ has at most $\aleph_1$ many equivalence classes.
\end{fact}

\begin{proof}
See \cite{Descriptive-Set-Theory-and-Infinitary-Languages}.
\end{proof}

The above fact may suggest that the properness of $\bbP_I$ should be used with countable support iterations to change covering numbers. It will be shown below that descriptive set theoretic techniques will give a stronger result in just $\mathsf{ZFC}$. However, in the context of proper forcing, the following combinatorial approach is worth mentioning:

\Begin{definition}{cov star of ideal}
Let $I$ be a $\sigma$-ideal on a Polish space $X$. $\text{cov}^*(I)$ is the smallest cardinal $\kappa$ such that there exists some $I^+$ $\bDelta_1^1$ $B \subseteq X$ and a set $U \subseteq I$ with $|U| = \kappa$ and $B \subseteq \bigcup U$. 
\end{definition}

\Begin{proposition}{cov greater than omega1 canonicalization borel}
Suppose $I$ is a $\sigma$-ideal such that $\text{cov}^*(I) > \omega_1$. Then $\bSigma_1^{1 \ \text{thin} \ \bDelta_1^1} \rightarrow_I \{\text{ev}\}$.
\end{proposition}

\begin{proof}
Let $\{C_\alpha : \alpha < \omega_1\}$ enumerate all the $E$-classes in order type $\omega_1$, using Fact \ref{number of classes of thin analytic equivalence relation}. Each $C_\alpha$ is $\bDelta_1^1$. 

Let $E \in \bSigma_1^{1 \ \text{thin} \ \bDelta_1^1}$. Let $B$ be an arbitrary $I^+$ $\bDelta_1^1$ set. $B = \bigcup_{\alpha < \omega_1} B \cap C_\alpha$. Since $\text{cov}^*(I) > \aleph_1$, there is some $\alpha$ such that $B \cap C_\alpha$ is $I^+$ $\bDelta_1^1$. Then $E \upharpoonright (B \cap C_\alpha) = \text{ev} \upharpoonright B \cap C_\alpha$.
\end{proof}

\Begin{proposition}{PFA thin analytic canonicalized borel}
If $\mathsf{PFA}$ holds, then for all $I$ such that $\bbP_I$ is proper, $\Sigma_1^{1 \ \text{thin} \ \bDelta_1^1} \rightarrow_I \{\text{ev}\}$. 
\end{proposition}

\begin{proof}
Let $B$ be a $I^+$ $\bDelta_1^1$ set. Let $U = \{C_\beta : \beta < \omega_1\}$ be a collection of $\bDelta_1^1$ sets in $I$. $\bbP_I$ being proper implies that $\bbP_I \upharpoonright B$ is proper. Let $D_\beta := \{F \in \bbP_I  \upharpoonright B : F \cap C_\beta = \emptyset\}$. $D_\beta$ is dense in $\bbP_I \upharpoonright B$. By $\mathsf{PFA}$, there is a filter $G \subseteq \bbP_I \upharpoonright B$ which is generic for $\{D_\beta : \beta < \omega_1\}$. $H$ constructs a real $x_H \in B$. By genericity, $x_H \notin C_\beta$ for all $\beta < \omega_1$. So $U$ can not cover $B$. $\text{cov}^*(I) > \aleph_1$. The result follows from Proposition \ref{cov greater than omega1 canonicalization borel}.
\end{proof}

The results are unsatisfactory in several ways. Models of $\mathsf{PFA}$ satisfy $\neg\mathsf{CH}$ and this was an essential fact since the proof used $\text{cov}^*(I) > \aleph_1$. Definability of the equivalence relation was not used in any deep way. The core of the proofs was combinatorial, using $\text{cov}^*(I) > \aleph_1$. 

The rest of this section provides results addressing the main question for thin $\bSigma_1^1$ equivalence relations which rely on definability properities of these equivalence relations. The best validation of the definability approach to thin $\analytic$ equivalence is that a stronger result will be proved with weaker assumptions (just $\mathsf{ZFC}$).

\Begin{fact}{thin equiv all name pinned under dense}
Let $E$ be a thin $\analytic$ equivalence relation on a Polish space $X$. Let $\bbP$ be a forcing. Suppose $\tau \in V^\bbP$ is such that $1_\bbP \forces_\bbP \tau \in X$. Then there is a dense set $D_\tau^E$ such that for all $p \in D_\tau^E$, $(p,p) \forces_{\bbP \times \bbP} \tau_\text{left} \ E \ \tau_\text{right}$, where $\tau_\text{left}$ and $\tau_{\text{right}}$ denote the $\bbP\times \bbP$ name for the evaluation of $\tau$ according to the left and right generic for $\bbP$, respectively, coming from a generic for $\bbP \times \bbP$. 
\end{fact}

\begin{proof}
This is due to Silver. See \cite{Effective-Enumeration-of-Classes-in-a}, Lemma 2.1 or the proof of \cite{Thin-Equivalence-Relations-and-Effective-Decompositions}, Theorem 2.3. 

A sketch of the result is provided:

Suppose not. Then there exists some $u \in \bbP$ such that for all $q \leq_\bbP u$, $(q,q) \not\forces_{\bbP \times \bbP} \tau_\text{left} \ E \ \tau_\text{right}$. Hence, there is some $u$ such that for all $q \leq_\bbP u$, there exists $q_0, q_1 \leq_\bbP q$ with $(q_0, q_1) \forces_{\bbP \times \bbP} \neg(\tau_\text{left} \ E \ \tau_\text{right})$. 

Suppose $E$ is a thin $\lanalytic(z)$ equivalence relation. Let $\Theta$ be some large ordinal such that $V_\Theta$ reflects the necessary statements to perform the proof below:

Let $N \prec V_\Theta$ be a countable elementary substructure with $z, \bbP, u, \tau \in N$. Let $M$ be the Mostowski collapse of $N$ with $\pi : N \rightarrow M$ be the Mostowski collapsing map. One may assume that for all $x$, $\text{tc}(x) \subseteq \omega$. So $\pi$ does not move reals or elements of $X$. In particular $\pi(z) = z$. Let $\bbQ = \pi(\bbP)$, $v = \pi(u)$, and $\sigma = \pi(\tau)$. By elementarity, $M$ satisfies that for all $q \leq_\bbQ v$, there exists $q_0, q_1 \leq_{\bbQ} q$ such that $(q_0, q_1) \forces_{\bbQ \times \bbQ} \neg(\sigma_\text{left} \ E \ \sigma_\text{right})$. 

Let $(D_n : n \in \omega)$ enumerate all the dense open sets in $\bbQ \times \bbQ$ of $M$. One may assume that $D_{n + 1} \subseteq D_n$, by replacing $D_n$ with $E_n = \bigcap_{m \leq n} D_m$. Next, a function $f : \finBinarySequence \rightarrow \bbQ$ will be constructed with the following properties:

(1) If $s \subseteq t$, then $f(t) \leq_\bbQ f(s)$. 

(2) For all $n \in \omega$, if $|s| = |t| = n$ and $s \neq t$, then $(f(s), f(t)) \in D_n$.

(3) For all $s \in \finBinarySequence$, $(f(s0), f(s1)) \forces_{\bbQ \times \bbQ} \neg(\sigma_\text{left} \ E \ \sigma_\text{right})$. 

\noindent To construct this $f$: Let $f(\emptyset) = v$. 

Suppose for all $s \in {}^n2$, $f(s)$ has been constructed with the above properties. For each $s \in {}^n2$, find some $q^{s0}, q^{s1} \leq (f(s), f(s))$ such that $(q^{s0}, q^{s1}) \forces_{\bbQ \times \bbQ} \neg(\sigma_\text{left} \ E \ \sigma_\text{right})$. Using the fact that $D_{n + 1}$ is dense open, for each $t \in {}^{n + 1}2$, find some $r^t \leq_\bbQ q^t$ such that for all $a,b \in {}^{n + 1}2$ with $a \neq b$, $(r^a, r^b) \in D_{n + 1}$. For $t \in {}^{n + 1}2$, define $f(t) = r^t$. 

For each $x \in \cantorspace$, let $G_x := \{p \in \bbQ : (\exists n)(f(x \upharpoonright n) \leq_\bbQ p)\}$. If $x,y \in \cantorspace$ and $x \neq y$, then $G_x \times G_y$ is $\bbQ \times \bbQ$-generic over $M$, using (2) and the assumption that for all $n \in \omega$, $D_{n + 1} \subseteq D_n$. So let $n$ be largest such that $x \upharpoonright n = y \upharpoonright n$. Let $s = x \upharpoonright n$. Without loss of generality, suppose $x(n) = 0$ and $y(n) = 1$. Then $f(s0) \in G_x$ and $f(s1) \in G_y$. Also $(f(s0), f(s1)) \forces_{\bbQ \times \bbQ} \neg(\sigma_\text{left} \ E \ \sigma_\text{right})$. By the forcing theorem applied in $M$, $M[G_x][G_y] \models \neg(\sigma[G_x] \ E\ \sigma[G_y])$. By Mostowski absoluteness, $V \models \neg(\sigma[G_x] \ E \ \sigma[G_y])$. 

Define $\Phi : \cantorspace \rightarrow X$ by $\Phi(x) = \sigma[G_x]$. By an appropriate coding, $\Phi$ is a $\borel$ function. $\Phi[\cantorspace]$ is a $\analytic$ set of pairwise disjoint $E$-inequivalent elements. By the perfect set property for $\analytic$ sets, there is a perfect set of pairwise $E$-inequivalent elements. This contradicts $E$ being a thin equivalence relation.
\end{proof}

\Begin{fact}{generic reals below condition are related}
Let $E$ be a thin $\analytic$ equivalence relation on a Polish space $X$. Let $\bbP$ be some forcing and $\tau \in V^\bbP$ be such that $1_{\bbP} \forces_\bbP \tau \in X$. Suppose $p \in D_\tau^E$. Let $M \prec H_\Theta$ be a countable elementary substructure with $\Theta$ sufficiently large and $\bbP, p, \tau \in M$. Then for all $G,H \in V$ such that $p \in G$, $p \in H$, and $G$ and $H$ are $\bbP$-generic over $M$, $V \models \tau[G] \ E \ \tau[H]$. 
\end{fact}

\begin{proof}
This is due to Silver. See \cite{Effective-Enumeration-of-Classes-in-a}, Lemma 2.4.

Suppose $G$ and $H$ are any two such generics. Let $K$ be such that it is $\bbP$-generic over $M[G][H]$. Then $M[G][K] \models \tau[G] \ E \ \tau[K]$ and $M[H][K] \models \tau[H] \ E \ \tau[K]$. By Mostowski absoluteness, $M[G][H][K] \models \tau[G] \ E \ \tau[H]$. By Mostowski absoluteness again, $V \models \tau[G] \ E \ \tau[H]$. 
\end{proof}

\Begin{theorem}{thin analytic canonicalize ev}
$\bSigma_1^{1 \ \text{thin}} \rightarrow_I \{\text{ev}\}$ whenever $I$ is a $\sigma$-ideal such that $\bbP_I$ is proper.
\end{theorem}

\begin{proof}
Let $B \in \bbP_I$. Since $D_{\dot x_\text{gen}}^E$ is dense, there exists some $B' \leq_{\bbP_I} B$ such that $B' \in D_{\dot x_\text{gen}}^E$. So $(B', B') \forces_{\bbP_I \times \bbP_I} (\dot x_\text{gen}) \ E \ (\dot x_\text{gen})$. Let $M \prec H_\Theta$ with $\Theta$ sufficiently large and $\bbP, B' \in M$. By Fact \ref{properness equivalence}, the set $C \subseteq B'$ of $\bbP_I$-generic over $M$ reals is a $I^+$ $\borel$ set. For $x \in C$, if $G_x$ denotes the $\bbP_I$-generic over $M$ filter constructed from $x$, then $\dot x_\text{gen}[G_x] = x$. Note that for all $x \in C$, $B' \in G_x$. By Fact \ref{generic reals below condition are related}, for all $x,y \in C$, $V \models x \ E \ y$. Hence $E \upharpoonright C = \text{ev} \upharpoonright C$. 
\end{proof}

Note that in this result, $E$ does not need to be an equivalence relation with all $\borel$ classes.

\section{Positive Answer for $\bPi_1^1$ Equivalence Relations}\label{Consistency of Positive Answer for bPi11 Equivalence Relations}
The variant of Question \ref{main question} for $\bPi_1^1$ equivalence relations can be phrased as follows:

\Begin{question}{main question pi11}
Let $\bPi_1^{1 \ \bDelta_1^1}$ be the class of $\bPi_1^1$ equivalence relations with all classes $\bDelta_1^1$ defined on $\bDelta_1^1$ subsets of Polish spaces. If $I$ is a $\sigma$-ideal on a Polish space $X$ such that $\bbP_I$ is proper, then does $\bPi_1^{1 \ \bDelta_1^1} \rightarrow_I \bDelta_1^1$ hold?
\end{question}

A positive answer for the $\coanalytic$ case follows from the same assumptions as the main question for $\analytic$ in almost the exact same manner as above:

\Begin{lemma}{pi11 increasing union borel equiv}
Let $E$ be a $\Pi_1^1(z)$ equivalence relation on $\bairespace$. Then there exists $\bDelta_1^1$ relations $E_\alpha$, for $\alpha < \omega_1$, with the property that if $\alpha < \beta$, then $E_\alpha \subseteq E_\beta$, $E = \bigcup_{\alpha < \omega_1} E_\alpha$ and there exists a club set $C \subseteq \omega_1$ such that for all $\alpha \in C$, $E_\alpha$ is an equivalence relation.
\end{lemma}

\begin{proof}
Let $T$ be a $z$-recursive tree on $\omega \times \omega \times \omega$ such that $(x,y) \in E \Leftrightarrow T^{(x,y)}$ is wellfounded. For each $\alpha < \omega_1$, define $E_\alpha := \{(x,y) : \text{rk}(T^{(x,y)}) < \alpha\}$. Each $E_\alpha$ is $\bDelta_1^1$. If $\alpha < \beta$, $E_\alpha \subseteq E_\beta$. $E = \bigcup_{\alpha < \omega_1} E_\alpha = E$. 

Let $C$ be the set of all $\alpha$ such that $E_\alpha$ is an equivalence relation. Increasing union of equivalence relations are equivalence relations so $C$ is closed. Fix $\alpha < \omega_1$. The set $D = \{(x,x) : x \in \bairespace\}$ is $\Sigma_1^1$. So by the boundedness theorem, there exist some $\delta < \omega_1$ such that $\text{rk}(T^{(x,x)}) < \delta$ for all $x \in \bairespace$. Let $\beta_0 = \max\{\alpha, \delta\}$. Suppose $\beta_n$ has been defined. The set $G = \{(x,y) : (y,x) \in E_{\beta_n}\}$ is $\bSigma_1^1$. By the boundedness theorem, there exists some $\beta' > \beta_n$ such that for all $(x,y) \in G$, $\text{rk}(T^{(x,y)}) < \beta'$. The set $H = \{(x,z) : (\exists y)((x,y) \in E_{\beta_n} \wedge (y, z) \in E_{\beta_n}\}$ is $\bSigma_1^1$. Again by the boundedness theorem, there exists some $\beta_{n + 1} > \beta'$ such that for all $(x,z) \in H$, $\text{rk}(T^{(x,z)}) < \beta_{n + 1}$. One has constructed an increasing sequence $\{\beta_n : n \in \omega\}$. Let $\beta = \sup\{\beta_n : n \in \omega\}$. Then $E_\beta$ is an equivalence relation. $C$ is unbounded.
\end{proof}

\Begin{lemma}{coanalytic equivalence stabilization of classes}
Let $E$ be a $\Pi_1^1(z)$equivalence relation. Let $x,y \in \bairespace$ be such that $[x]_E$ is $\Sigma_1^1(y)$. Let $\delta$ be an ordinal such that $\omega_1^{y \oplus z} \leq \delta$ and $E_\delta$ is a equivalence relation. Then $[x]_E = [x]_{E_\delta}$. 
\end{lemma}

\begin{proof}
Define $E' \subseteq (\bairespace)^2$ by
$$a \ E' \ b \Leftrightarrow (a \in [x]_E \wedge b \in [x]_E) \vee (a = b)$$
$E'$ is $\Sigma_1^1(y)$. $E' \subseteq E$. By the effective boundedness theorem, there exists a $\alpha < \omega_1^{y \oplus z} \leq \delta$ such that for all $(x,y) \in E'$, $\text{rk}(T^{(x,y)}) < \alpha$. Hence $E' \subseteq E_\alpha$. 

Since $E_\delta \subseteq E$, $[x]_\delta \subseteq [x]_E$. Also $[x]_E = [x]_{E'} \subseteq [x]_{E_\alpha} \subseteq [x]_{E_\delta}$. Therefore, $[x]_E = [x]_{E_\delta}$. 
\end{proof}

Note that the parameter $x$ is not used in the above lemma. This is in contrast to Lemma \ref{analytic equivalence stabilization of classes}. Surprisingly, this observation will be used later. (See Proposition \ref{speculation counterexample proposition}.)

\Begin{lemma}{true sigma11 code in generic extensions}
Let $E$ be a $\Pi_1^1(z)$ equivalence relation. Let $M \prec H_\Theta$ be a countable elementary substructure with $z \in M$. Let $\bbP$ be a forcing in $M$ which adds a generic real. Suppose for all $g$ which are $\bbP$-generic over $M$, there exists a $y \in M[g]$ such that $V \models [g]_E$ is $\Sigma_1^1(y)$. Then there exists a countable ordinal $\alpha$ such that for all $\bbP$-generic over $M$ reals $g$, $[g]_E = [g]_{E_\alpha}$. 
\end{lemma}

\begin{proof}
The proof is almost the same as the proof for Lemma \ref{true sigma11 code in generic extensions} using Lemma \ref{coanalytic equivalence stabilization of classes} in place of Lemma \ref{analytic equivalence stabilization of classes}.
\end{proof}

These previous results can be used to give a positive answer for a specific $\Pi_1^1$ equivalence relation in $\mathsf{ZFC}$.

\Begin{example}{hyp arrow borel}
Let $H$ be an equivalence relation on $\bairespace$ defined by $x \ H \ y$ if and only if $x \in L_{\omega_1^y}(y) \wedge y \in L_{\omega_1^x}(x)$. 

$H$ is a $\Pi_1^1$ equivalence relation with all classes countable. $H$ is the equivalence relation of hyperarithmetic equivalence.

If $I$ is a $\sigma$-ideal on $\bairespace$ with $\bbP_I$ proper, then $\{H\} \rightarrow_I \bDelta_1^1$. 
\end{example}

\begin{proof}
Fix $B$ an $I^+$ $\bDelta_1^1$ set. Choose $M \prec H_\Theta$ with $\Theta$ sufficiently large and $B, \bbP_I \in M$. By Fact \ref{properness equivalence}, let $C \subseteq B$ be the set of $\bbP_I$-generic over $M$ elements in $B$. Let $x \in C$. $\omega_1^x$ is a countable ordinal in $M[x]$. In $M[x]$, $L_{\omega_1^x}(x)$ is countable. $[x]_E^{M[x]} \subseteq L_{\omega_1^x}(x)$. Therefore, in $M[x]$, there is a function $f : \omega \rightarrow \bairespace$ such that $f$ enumerates $[x]_H^{M[x]}$. By absoluteness, $[x]_H = [x]_H^{M[x]}$. So $[x]_H$ is $\Delta_1^1(f)$ and $f \in M[x]$. By Lemma \ref{true sigma11 code in generic extensions}, there is some countable ordinal $\alpha$ such that $[x]_E = [x]_{E_\alpha}$ for all $x \in C$. So $E \upharpoonright C = E_\alpha \upharpoonright C$. $E \upharpoonright C$ is $\bDelta_1^1$. 
\end{proof}

\Begin{definition}{set of codes for sigma11 sets}
Let $E$ be a $\Pi_1^1(z)$ equivalence relation. Define the set $D$ by
$$(x,T) \in D \Leftrightarrow (T \text{ is a tree on } \omega \times \omega) \wedge (\forall y)(y \ E \ x \Leftrightarrow T^y \notin \text{WF})$$
$D$ is $\Pi_2^1(z)$.
\end{definition}

\Begin{theorem}{universally baire pi14 generic absoluteness coanalytic arrow borel}
Assume all $\bPi_2^1$ sets are universally Baire and $\bPi_4^1$-generic absoluteness holds. Let $I$ be a $\sigma$-ideal such that $\bbP_I$ is proper. Then $\bPi_1^{1 \ \bDelta_1^1} \rightarrow_I \bDelta_1^1$. 
\end{theorem}

\begin{proof}
The proof is the same as Theorem \ref{universally baire pi14 generic absoluteness analytic arrow borel} with the required change.
\end{proof}

A similar argument using iterable structures as in the $\analytic$ case yields a positive answer from a more precise assumption with lower consistency strength. 

\Begin{proposition}{all classes borel coanalytic case}
Let $E$ be a $\Pi_1^1(z)$ equivalence relation. There is a $\Pi_3^1(z)$ formula $\varpi(v)$ in free variable $v$ such that:

Let $x \in \bairespace$. If $(x \oplus z)^\sharp$ exists, then the statement ``$[x]_E$ is $\borel$'' is equivalent to $\varpi(x)$. 

Assume for all $r \in \bairespace$, $r^\sharp$ exists. The statement ``all $E$-classes are $\borel$'' is equivalent to $(\forall x)\varpi(x)$. In particular, this statement is $\Pi_3^1(z)$. 
\end{proposition}

\begin{proof}
Assume for simplicity, $E$ is a $\lcoanalytic$ equivalence relation on $\bairespace$. Let $T$ be a tree on $\omega \times \omega \times \omega$ such that 
$$(x,y) \in E \Leftrightarrow T^{x,y} \text{ is wellfounded}$$
Let $\varpi(v)$ be the statement:
$$(\forall y)(y = v^\sharp \Rightarrow ``1_{\text{Coll}(\omega, < c_1)} \forces_{\text{Coll}(\omega, < c_1)} (\exists c)(c \in \text{WO} \wedge (\forall y)((y \ E \ v) \Rightarrow \text{rk}(T^{v,y}) < \text{ot}(c)))" \in y)$$
The rest of the argument is the same as in Proposition \ref{sharps of reals borelness of all classes pi13}.
\end{proof}

\Begin{definition}{coanalytic name associated with all classes borel}
Let $I$ be a $\sigma$-ideal on a Polish space $X$ such that $\bbP_I$ is proper. Let $\mu_E^I$ be $\tau_{\neg \varpi}^{\bbP_I}$ from Fact \ref{name associated statement}.
\end{definition}

\Begin{definition}{coanalytic name associated all classes pi11}
Let $I$ be a $\sigma$-ideal on $\bairespace$ such that $\bbP_I$ is proper.  Consider the formula ``$(\exists y)([\dot x_\text{gen}]_E \text{ is } \Sigma_1^1(y))$''. Write it as $(\exists y)\psi(y)$. By Fact \ref{name associated statement}, let $\sigma_E^{I}$ be $\tau_\psi^{\bbP_I}$. 
\end{definition}

\Begin{definition}{associated triple for consistency main question coanalytic}
Suppose $I$ is a $\sigma$-ideal on $\bairespace$ such that $\bbP_I$ is proper. Let $E \in \coanalytic^\borel$. Define $\chi_E^I = \langle \bbP_I, \mu_E^{\bbP_I}, \sigma_E^{\bbP_I}\rangle$. 
\end{definition}

\Begin{theorem}{sharps of reals and positive answer main question coanalytic}
Suppose $I$ is a $\sigma$ ideal on $\bairespace$ such that $\bbP_I$ is proper. If for all $x \in \cantorspace$, $x^\sharp$ exists and $(\chi_E^I)^\sharp$ exists for all $E \in \coanalytic^\borel$, then $\coanalytic^\borel \rightarrow_I \borel$. 
\end{theorem}

\begin{proof}
This is similar to Theorem \ref{sharps of reals and positive answer main question}.
\end{proof}

\Begin{corollary}{sharp hereditarily continuum plus canonicalization coanalytic}
If $z^\sharp$ exists for all $z \in H_{(2^{\aleph_0})^+}$, then $\coanalytic^\borel \rightarrow_I \borel$ for all $\sigma$-ideal $I$ such that $\bbP_I$ is proper.
\end{corollary}

\section{$\bPi_1^1$ Equivalence Relations with Thin or Countable Classes}\label{bPi11 Equivalence Relations with Thin or Countable Classes}

The preservation of the statement ``all classes are $\bDelta_1^1$'' played an important role in the consistency results above. Next, one will consider $\bPi_1^1$ equivalence relations which are very sensitive to set theoretic assumptions and generic extensions.

\Begin{definition}{thin sets}
Let $X$ be a Polish space. $A \subseteq X$ is thin if and only if it does not contain a perfect set.
\end{definition}

\Begin{fact}{largest thin pi11 set}
For each $z \in \cantorspace$, define $Q_z := \{x \in \cantorspace : x \in L_{\omega_1^{z \oplus x}}(z)\}$. $Q_z$ is the largest thin $\Pi_1^1(z)$ set in the sense that if $S$ is a thin $\Pi_1^1(z)$ set, then $S \subseteq Q_z$. Moreover, for each $\alpha < \omega_1^{L[z]}$, there exists some $x \in Q_z$ such that $\alpha < \omega_1^x$. Therefore, if $\omega_1^{L[z]} = \omega_1$, then $Q_z$ is an uncountable thin $\Pi_1^1(z)$ set. It is consistent that $\bPi_1^1$ sets do not have the perfect set property.
\end{fact}

\begin{proof}
See \cite{Recursive-Aspects-of-Descriptive-Set-Theory}, pages 83-87. One will give the $\Pi_1^1(z)$ definition to get a better understanding of what $Q_z$ is:
$$x \in Q_z \Leftrightarrow (\forall M)((M \text{ is an $\omega$-model of }\mathsf{KP} \wedge z \in M \wedge x \in M) \Rightarrow M \models x \in L[z])$$
So $x \in Q_z$ if and only if $L_{\omega_1^{z \oplus x}}(z \oplus x) = L_{\omega_1^{z \oplus x}}(z)$. Or put another way, the smallest admissible set containing $z \oplus x$ is a model of $V = L[z]$. Certainly $Q_z \subseteq L[z]$. So $Q_z$ can also be thought of as the set of reals that appear in $L[z]$ very quickly in the sense that $x \in Q_z$ if and only if the first ordinal $\alpha$ such that $L_{\alpha}(z \oplus x)$ is admissible is also the first $z$-admissible ordinal $\alpha$ such that $x \in L_\alpha[z]$. 
\end{proof}

Now one can give a simple example of an equivalence relation $E$, a model of $\mathsf{ZFC}$, and forcing which does not preserve the statement ``all $E$ classes are $\bDelta_1^1$''. Note that this statement is $\bPi_4^1$ so it will be preserved if the universe satisfies $\bPi_4^1$-generic absoluteness. The desired example will necessarily have to reside in a universe with weak large cardinals.

\Begin{definition}{inaccessible to the reals}
$\omega_1$ is inaccessible to reals if and only if for all $x \in \bairespace$, $L[x] \models (\omega_1^V \text{ is inaccessible})$ if and only if for all $x \in \bairespace$, $\omega_1^{L[x]} < \omega_1$.
\end{definition}

\Begin{proposition}{all classes borel not preserved}
Let $E$ and $F$ be equivalence relations on $(\bairespace)^2$ defined by
$$(a,x) \ E \ (b, y) \Leftrightarrow (a = b) \wedge (x, y \in Q_a \vee x = y)$$
$$(a,x) \ F \ (b,y) \Leftrightarrow (a = b) \wedge (x,y \notin Q_a \vee x = y)$$
$E$ is a $\Pi_1^1$ equivalence relation and $F$ is $\Sigma_1^1$ equivalence relation. 

Let $\kappa$ be an inaccessible cardinal in $L$. Let $G \subseteq \text{Coll}(\omega, <\kappa)$ be $\text{Coll}(\omega, <\kappa)$-generic over $L$. Then $L[G] \models$ $E$ and $F$ have all classes $\bDelta_1^1$. Let $g \subseteq \text{Coll}(\omega, \kappa)$ be $\text{Coll}(\omega,\kappa)$-generic over $L[G]$, then $L[G][g] \models$ not all $E$ and $F$ classes are $\bDelta_1^1$.  
\end{proposition}

\begin{proof}
The formula provided in the proof of Fact \ref{largest thin pi11 set} shows that the formula ``$x \in Q_z$'' is $\Pi_1^1$ in variables $x$ and $z$. From this, it follows that $E$ and $F$ are $\Pi_1^1$ and $\Sigma_1^1$, respectively.

In $L[G]$, $\omega_1$ is inaccessible to reals (\cite{Set-Theory-Exploring} Theorem 8.20). For each $(a,b)$, $[(a,b)]_E$ is either a singleton or in bijection with $Q_a$. Since $Q_a \subseteq (\bairespace)^{L[a]}$, in all cases, $[(a,b)]_E$ is countable and hence $\bDelta_1^1$. $F$-classes are then singletons or complements of countable sets. All $F$-classes are $\bDelta_1^1$. 

$\text{Coll}(\omega, <\kappa) * \dot{\text{Coll}}(\omega, \check\kappa)$ is a forcing (in $L$) of size $\kappa$ which collapses $\kappa$ to $\omega$. Such forcing are forcing equivalent to $\text{Coll}(\omega, \kappa)$ by \cite{Set-Theory} Lemma 26.7. Let $h \subseteq \text{Coll}(\omega, \kappa)$ which is $\text{Coll}(\omega, \kappa)$-generic over $L$ with $L[h] = L[G][g]$. $L[G][g] \models \omega_1^{L[h]} = \omega_1$. $\omega_1$ is not inaccessible to reals in $L[G][g]$. Moreover, $[(h,h)]_E$ is not $\bDelta_1^1$ in $L[G][g]$ as it is an uncountable thin set and the perfect set property holds for $\bSigma_1^1$ sets. Simiarly, $F$ has a class which is not $\bDelta_1^1$. 
\end{proof}

In the previous example, in $L[G]$, $\text{Coll}(\omega, \omega_1^{L[G]}) = \text{Coll}(\omega, \kappa)$ is not a proper forcing. One may ask whether there is an $\analytic$ or $\coanalytic$ equivalence relation with all classes $\borel$ and a proper forcing coming from a $\sigma$-ideal on a Polish space such that in the induced generic extension, the statement that ``all classes are $\borel$'' is false. Sy-David Friedman's forcing to code subsets of $\omega_1$ is an $\aleph_1$-c.c. forcing which can be repesented as an idealize forcing which (like in the proof of the above proposition) adds a real $r$ such that $L[A][r] = L[r]$. The two equivalence relations from the above proposition can be used with this forcing to give a similar result. See Section \ref{Conclusion} for more details about this forcing.

For $\bSigma_1^1$ equivalence relation with all classes countable, Proposition \ref{all classes countable arrow borel} shows that the main question has a positive answer without additional set theoretic assumptions. There were two important aspects of the proof. First, the countability of all classes of a $\bSigma_1^1$ equivalence relation is $\bSigma_1^1$ and hence remains true in all generic extensions. This fact is used to give an enumeration $f$ of $[x]_E$ in $M[x]$. Secondly, the statement that $f$ enumerates $[x]_E$ is $\bPi_1^1$ and hence absolute between (the countable model) $M[x]$ and $V$.

One can ask the same question for $\bPi_1^1$ equivalence relations with all classes countable. However, the above proof can not be applied. First, the countability of all classes of a $\bPi_1^1$ equivalence relation is $\bPi_4^1$. Secondly, the statement that some function $f$ enumerates $[x]_E$ is $\bPi_2^1$; hence, it does not necessarily persist from $M[x]$ to $V$.

The $\bPi_1^1$ equivalence relations where these issues are most perceptible are the equivalence relations $E$ with all classes countable but for some $x$, $L[x] \models [x]_E \text{ is uncountable }$. It is not provable that all $E$-classes are countable; however, all the $E$-classes are thin.

\Begin{proposition}{pi11 equivalence relation all classes thin preserved}
Let $A$ be a $\bPi_1^1$ set. The statement ``$A$ is thin'' is $\bPi_2^1$. Let $E$ be a $\bPi_1^1$ equivalence relation. The statement ``all $E$-classes are thin'' is $\bPi_2^1$. Both of these statements are absolute to generic extensions.
\end{proposition}

\begin{proof}
$$(\forall T)(T \text{ is perfect tree } \Rightarrow ((\exists x)((\forall n)(x \upharpoonright n \in T) \wedge x \notin A)))$$
$$(\forall x)(\forall T)(T \text{ is perfect tree } \Rightarrow ((\exists y)((\forall n)(y \upharpoonright n \in T) \wedge \neg (x \ E \ y))))$$
These two $\bPi_2^1$ formulas are equivalent to ``$A$ is thin'' and ``all $E$ classes are thin'', respectively.
\end{proof}

\Begin{definition}{class coanalytic equiv all classes countable}
Let $\bPi_1^{1 \ \aleph_0}$ denote the class of all $\bPi_1^1$ equivalence relations with all classes countable defined on $\bDelta_1^1$ subsets of Polish spaces. Let $\bPi_1^{1 \ \text{thin}}$ denote the class of all $\bPi_1^1$ equivalence relation with all classes thin defined on $\bDelta_1^1$ subsets of Polish spaces.
\end{definition}

\Begin{theorem}{omega1 countable meager ideal coanalytic thin arrow borel}
If $\omega_1^L < \omega_1$, then $\bPi_1^{1 \ \text{thin}} \rightarrow_{I_\text{meager}} \bDelta_1^1$. 
\end{theorem}

\begin{proof}
Fix a non-meager $\bDelta_1^1$ set $B$. Let $\bbC$ denote Cohen forcing, i.e. finite partial functions from $\omega$ into $2$. Let $U$ be the set of all constructible dense subsets of $\bbC$. Since $\omega_1^L < \omega_1$, $|U| = \aleph_0$. 

Let $M \prec H_\Theta$ be a countable elementary substructure with $\Theta$ a sufficiently large cardinal, $B, \bbP_{I_\text{meager}}, \omega_1^L, U \in M$, $\omega_1^L \subseteq M$, and $U \subseteq M$. By Fact \ref{properness equivalence}, let $C$ be the set of $\bbP_{I_\text{meager}}$-generic over $M$ reals in $B$. 

Take $x \in C$. Since $\bbC$ and $\bbP_{I_\text{meager}}$ are forcing equivalent and $U \subseteq M$, $x$ is also $\bbC$-generic over $L$. Since $\bbC$ satisfies the $\aleph_1$-chain condition, $\omega_1^{L[x]} = \omega_1^L < \omega_1$. Since $\omega_1^L$ is countable in $M$, $L_{\omega_1^{L[x]}}[x] = L_{\omega_1^L}[x] \subseteq M[x]$ and is countable there. Since $[x]_E$ is thin, $[x]_E \subseteq L_{\omega_1^{L[x]}}[x]$. In $M[x]$, there is an enumeration $f : \omega \rightarrow ([x]_E)^{M[x]}$. The claim is that $([x]_E)^{M[x]} = [x]_E$: since $([x]_E)^{V} \subseteq L_{\omega_1^{L[x]}}[x] \subseteq M[x]$, $M[x] \models y \ E \ x \Leftrightarrow (L[x])^{M[x]} \models y \ E \ x \Leftrightarrow L[x] \models y \ E \ x \Leftrightarrow V \models y \ E \ x$, by Mostowski absoluteness.

 Therefore, in $V$, $[x]_E$ is $\Delta_1^1(f)$ and $f \in M[x]$. By Lemma \ref{true sigma11 code in generic extensions}, there is some countable $\alpha < \omega_1$, such that $E \upharpoonright C = E_\alpha \upharpoonright C$. The latter is $\bDelta_1^1$. 
\end{proof}

If $\omega_1$ is inaccessible to reals, then $\bPi_1^{1 \ \text{thin}} = \bPi_1^{1\ \aleph_0}$. Familiar models that satisfy $\omega_1$ is inaccessible to reals include generic extensions of the L\'{e}vy collapse of an inaccessible cardinal to $\omega_1$. Next, one will consider the main question for $\bPi_1^{1 \ \aleph_0}$ in models of this type and obtain some improved consistency results.

The main large cardinal useful here is the remarkable cardinal isolated in \cite{Proper-Forcing-and-Remarkable-Cardinals-II} to understand absoluteness for proper forcing in $L(\bbR)$. It is a fairly weak large cardinal. Its existence is consistent relative to $\omega$-Erd\"{o}s cardinals. Remarkable cardinals are compatible with $\mathsf{V = L}$. If $0^\sharp$ exists, then all Silver's indiscernibles are remarkable cardinals in $L$. Also if $\kappa$ is a remarkable cardinal, then $\kappa$ is a remarkable cardinal in $L$. 

\Begin{definition}{remarkable cardinal definition}
(\cite{Proper-Forcing-and-Remarkable-Cardinals-II} Definition 1.1) A cardinal $\kappa$ is a remarkable cardinal if and only if for all regular cardinals $\theta > \kappa$, there exists $M, N, \pi, \sigma, \bar\kappa$ and $\bar{\theta}$ such that the following holds:

(i) $M$ and $N$ are countable transitive sets.

(ii) $\pi : M \rightarrow H_\theta$ is an elementary embedding.

(iii) $\pi(\bar{\kappa}) = \kappa$

(iv) $\sigma : M \rightarrow N$ is an elementary embedding with $\text{crit}(\sigma) = \bar{\kappa}$. 

(v) $\bar\theta = \text{ON} \cap M$, $\sigma(\bar\kappa) > \bar\theta$, and $N \models \theta$ is a regular cardinal.

(vi) $M \in N$ and $N \models M = H_\theta$. 
\end{definition}

\Begin{fact}{levy collapse remarkable proper forcing real generic over L}
(Schindler) Let $\kappa$ be a remarkable cardinal in $L$. Let $G \subseteq \text{Coll}(\omega, < \kappa)$ be $\text{Coll}(\omega, <\kappa)$-generic over $L$. Let $\bbP \in L[G]$ be a proper forcing. Let $H \subseteq \bbP$ be a $\bbP$-generic filter over $L[G]$. If $x \in (\bairespace)^{L[G][H]}$, then there exists a forcing $\bbQ \in L_\kappa$ and a $K \subseteq \bbQ$ in $L[G][H]$ which is $\bbQ$-generic over $L$ and $x \in L[K]$. 
\end{fact}

\begin{proof}
See \cite{Proper-Forcing-and-Remarkable-Cardinals-II}, Lemma 2.1.
\end{proof}

\Begin{theorem}{remarkable cardinal consistency pi11aleph0 arrow borel}
Let $\kappa$ be a remarkable cardinal in $L$. Let $G \subseteq \text{Coll}(\omega, < \kappa)$ be $\text{Coll}(\omega, <\kappa)$ generic over $L$. In $L[G]$, if $I$ is a $\sigma$-ideal such that $\bbP_I$ is proper, then $\bPi_1^{1 \ \aleph_0} \rightarrow_I \bDelta_1^1$. 
\end{theorem}

\begin{proof}
Working in $L[G]$, let $E$ be an equivalence relation in $\bPi_1^{1 \ \aleph_0}$. For simplicity, let assume $E$ is $\Pi_1^1$ (otherwise one should include the parameter defining $E$ in all the discussions below). In particular, all $E$-classes are thin, and this statement will be absolute to all generic extensions.

Let $B$ be an $I^+$ $\bDelta_1^1$ set. Let $M \prec H_\Theta$ be a countable elementary substructure, $\Theta$ sufficiently large cardinal, and $B, \bbP_I, G \in M$. $H_\Theta = L_\Theta[G]$. Therefore, $M = L^M[G]$. Note that from the point of view of $M$, $G$ is $L^M$ generic for $\text{Coll}(\omega, <\kappa)^M$. Using Fact \ref{properness equivalence}, let $C \subseteq B$ be the $I^+$ $\bDelta_1^1$ set of $\bbP_I$-generic over $M$ reals in $B$.

Fix $x \in C$. Applying Fact \ref{levy collapse remarkable proper forcing real generic over L} in $M[x] = L^M[G][x]$, there exists some $\bbQ \in (L_\kappa)^M$ and $K \subseteq \bbQ$ in $M[x]$ which is $\bbQ$-generic over $L^M$ such that $x \in L[K]$. Since $M$ satisfies $\bbQ \in L_\kappa$ and $\kappa$ is a remarkable cardinal (in particular inaccessible) in $L$, $M$ thinks that $\mathscr{P}^L(\bbQ) \in L_\kappa$. Since $M = L[G]$, $M \models \mathscr{P}^L(\bbQ)$ is countable. Let $f : \omega \rightarrow \mathscr{P}^L(\bbQ)$ be a function in $M$ such that $M$ thinks it surjects onto $\mathscr{P}^L(\bbQ)$. Since $M \prec (H_\Theta)^{L[G]}$, $f$ really is a surjection of $\mathscr{P}^L(\bbQ)$ in the real universe $L[G]$. This establishes that $\mathscr{P}^L(\bbQ) \subseteq M$. In particular, $\mathscr{P}^L(\bbQ) \subseteq L^M$. This and the fact that $K$ is generic over $L^M$ imply that $K$ is $\bbQ$-generic over the real $L$. Since $\bbQ \in L_\kappa$, all cardinals of $L$ greater than $|\bbQ|$ is preserved in $L[K]$. Therefore, $\omega_1^{L[x]} \leq \omega_1^{L[K]} \leq (|\bbQ|^+)^L$. Since $\bbQ \in M$ and $M \models \bbQ \in L$, there is some ordinal $\alpha$ such that $M \models L \models |\bbQ|^+ = \alpha$. Because $M \prec H_\Theta$, the real universe $L[G]$ satisfies $L \models |\bbQ|^+ = \alpha$. This establishes that $(|\bbQ|^+)^L \in M$. Since $\bbQ \in L_\kappa$ and $\kappa$ is inaccessible, $(|\bbQ|^+)^L < \kappa$. Since $M = L^M[G]$, $(|\bbQ|^+)^L$ is a countable ordinal in $M$. As shown above, $\omega_1^{L[x]} < (|\bbQ|^+)^L$, so in $M$, $\omega_1^{L[x]}$ is countable. Since $[x]_E$ is thin, $[x]_E \subseteq (\bairespace)^{L[x]}$. As $\omega_1^{L[x]}$ is countable in $M[x]$, $M[x] \models [x]_E$ is countable. There exists some surjection $h : \omega \rightarrow ([x]_E)^{M[x]}$. The claim is that $([x]_E)^{L[G]} = ([x]_E)^{M[x]}$: since $([x]_E)^{L[G]} \subseteq L_{\omega_1^{L[x]}}[x] \subseteq M[x]$, $M[x] \models y \ E \ x \Leftrightarrow (L[x])^{M[x]} \models y \ E \ x \Leftrightarrow (L[x])^{L[G]} \models y \ E \ x \Leftrightarrow L[G] \models y \ E \ x$, by Mostowski absoluteness. 

Therefore, in $L[G]$, $[x]_E$ is $\Delta_1^1(h)$ and $h \in M[x]$. By Lemma \ref{true sigma11 code in generic extensions}, there is a countable $\alpha < \omega_1^{L[G]}$, such that $E \upharpoonright C = E_\alpha \upharpoonright C$. $E_\alpha \upharpoonright C$ is $\bDelta_1^1$. 
\end{proof}

Using some well-known results of Kunen, a similar proof shows that the consistency of $\bPi_1^{1 \ \aleph_0} \rightarrow_I \bDelta_1^1$ for $I$ such that $\bbP_I$ is $\aleph_1$-c.c. follows from the consistency of a weakly compact cardinal.

\Begin{fact}{weakly compact element of extension generic for small forcing}
(Kunen) Let $\kappa$ be a weakly compact cardinal. Let $\bbP$ be a $\kappa$-c.c. forcing. Let $G \subseteq \bbP$ be $\bbP$-generic over $V$. If $x \in H_\kappa^{V[G]}$, then there exists a forcing $\bbQ \in V_\kappa$ and a $K \subseteq \bbQ$ which is generic over $V$ such that $x \in V[K]$. 
\end{fact}

\begin{proof}
This is due to Kunen. See \cite{The-Necessary-Maximality-Principle-for-C.C.C.}, Lemma 5.3 for a proof.
\end{proof}

\Begin{theorem}{weakly compact consistency pi11aleph0 arrow borel for c.c.c.}
Let $\kappa$ be a weakly compact cardinal in $L$. Let $G \subseteq \text{Coll}(\omega, <\kappa)$ be $\text{Coll}(\omega, <\kappa)$-generic over $L$. In $L[G]$, if $I$ is a $\sigma$-ideal such that $\bbP_I$ is $\aleph_1$-c.c., then $\bPi_1^{1 \ \aleph_0} \rightarrow_I \bDelta_1^1$. 
\end{theorem}

\begin{proof}
Let $\dot \bbP_I$ be a name for $\bbP_I$ in $L[G]$. $\text{Coll}(\omega, <\kappa)$ satisfies the $\kappa$-chain condition. Since $\aleph_1^{L[G]} = \kappa$, for some $p \in G$, $p \forces_{\text{Coll}(\omega, <\kappa)} \dot \bbP_I$ satisfies the $\check \kappa$-chain condition. By considering the forcing of conditions below $p$, one may as well assume $p = 1_{\text{Coll}(\omega, < \kappa)}$. Then $\text{Coll}(\omega, < \kappa) * \dot \bbP_I$ satisfies the $\kappa$-chain condition. Now use Fact \ref{weakly compact element of extension generic for small forcing} and finish the proof much like in Theorem \ref{remarkable cardinal consistency pi11aleph0 arrow borel}.
\end{proof}

\section{$\bDelta_2^1$ Equivalence Relations with all Classes $\borel$}\label{bPi21 Equivalence Relations with all bDelta11 Classes}

One can ask the same question for $\bDelta_2^1$ equivalence relation with all classes $\borel$: If $E$ is a $\bDelta_2^1$ equivalence relation with all classes $\borel$ and $I$ is a $\sigma$-ideal such that $\bbP_I$ is proper, does $\{E\} \rightarrow_I \borel$ hold?

It will be shown that in $L$, there is a $\bDelta_2^1$ equivalence relation $E_L$ such that $\{E_L\} \rightarrow_I \bDelta_1^1$ does not hold for any $\sigma$-ideal $I$. 

\Begin{definition}{E L equivalence relation}
Let $E_L$ be the equivalence relation defined on $\cantorspace$ by $x \ E_L \ y$ if and only if
$$(\forall A)((A \text{ is a well-founded $\omega$-model of } \mathsf{KP + V = L}) \Rightarrow (x \in A \Leftrightarrow y \in A))$$
$E_L$ is a $\Pi_2^1$ equivalence relation.
\end{definition}

Note that $A$ is a structure with domain $\omega$. As $A$ is an $\omega$-model, there is an isomorphic copy of $\omega$ in $A$. The statement ``$x \in A$'' should be understood using this copy of $\omega$ in $A$.

Rather than $\mathsf{KP + V = L}$, one could also use some $\Upsilon + \mathsf{V = L}$ where $\Upsilon$ is a large enough fragment to $\mathsf{ZFC}$ to perform the forcing argument below. If one is willing to assume that there exists a transitive model of $\mathsf{ZFC}$, then one can replace the above with $\mathsf{ZFC + V = L}$ and be in the familiar setting.

\Begin{definition}{admissible level real appear function}
Assume $V = L$. Let $\iota : \cantorspace \rightarrow \omega_1$ be the function such that $\iota(x)$ is the smallest admissible ordinal $\alpha$ such that $x \in L_\alpha$. 
\end{definition}

\Begin{proposition}{EL equivalent same iota}
For all $x, y \in \cantorspace$, $x \ E_L \ y$ if and only if $\iota(x) = \iota(y)$. 
\end{proposition}

\begin{proof}
Assume $\iota(x) = \iota(y)$. Let $A$ be a wellfounded model of $\mathsf{KP + V = L}$ such that $x \in A$. There is some $\beta$ such that $L_\beta$ is the Mostowski collapse of $A$. $L_\beta$ is transitive and satisfies $\mathsf{KP}$, so it is an admissible set. $\beta$ is an admissible ordinal. $\iota(x) \leq \beta$. $y \in L_{\iota(y)} = L_{\iota(x)} \subseteq L_\beta$. So $y \in A$. Hence $x \in A$ implies $y \in A$. By a symmetric argument, $y \in A$ implies $x \in A$. $x \ E_L \ y$. 

Assume $x \ E_L \ y$. Suppose $\alpha < \omega_1$ with $L_\alpha \models \mathsf{KP}$ and $x \in L_\alpha$. Since $L_\alpha$ is countable, there is a countable structure $A$ with domain $\omega$ isomorphic to $L_\alpha$. $A \models \mathsf{KP}$, $A$ is an $\omega$-model, and $x \in A$. $x \ E_L \ y$ implies that $y \in A$. Therefore, $y \in L_\alpha$. Hence $\iota(x) \leq \iota(y)$. By a symmetric argument, $\iota(y) \leq \iota(x)$. $\iota(x) = \iota(y)$. 
\end{proof}

Earlier drafts of this paper only asserted that $E_L$ was $\Pi_2^1$. Drucker observed that a very similar equivalence relation to $E_L$ was actually $\Delta_2^1$:

\Begin{proposition}{EL is delta12}
(Drucker) $E_L$ is $\Delta_2^1$. 
\end{proposition}

\begin{proof}
The claim is that
$$x \ E_L \ y \Leftrightarrow H_{\aleph_1} \models (\exists M)((M \text{ is transitive}) \wedge (x, y \in M) \wedge (M \models \mathsf{KP + V = L}) \wedge (M \models \psi(x,y)))$$
where 
$$\psi(x,y) \Leftrightarrow (\forall A)((A \text{ is transitive} \wedge A \models \mathsf{KP + V = L}) \Rightarrow (x \in A \Leftrightarrow y \in A))$$

To see this: $(\Rightarrow)$ By Proposition \ref{EL equivalent same iota}, $\iota(x) = \iota(y)$. Then $H_{\aleph_1}$ satisfies the above formula using $L_{\iota(x)}$. 

$(\Leftarrow)$ Suppose $\neg(x \ E_L \ y)$. Let $M$ witness the negation of the statement from Defintion \ref{E L equivalence relation}. Without loss of generality, $\iota(x) < \iota(y)$. By $\Delta_1$ absoluteness, if $H_{\aleph_1}$ thinks $M$ is transitive and satisfies $\mathsf{KP + V = L}$, then $M$ really is transitive and satisfies $\mathsf{KP + V + L}$. So $M = L_\alpha$ for some $\alpha < \omega_1$. Since $x,y \in M = L_{\alpha}$, $\alpha \geq \iota(y)$. Then $M \models \neg(\psi(x,y))$ since $L_{\iota(x)} \in L_{\alpha} = M$, $x \in L_{\iota(x)}$, and $y \notin L_{\iota(x)}$. 

$\psi(x,y)$ is a first order formula in the language of set theory. First order satisfaction is $\Delta_1$. The above shows that $x \ E_L \ y$ is equivalent to a formula which is $\Sigma_1$ over $H_{\aleph_1}$. Hence $E_L$ is $\Sigma_2^1$. 
\end{proof}

Assuming $\mathsf{V =L}$, Proposition \ref{EL equivalent same iota} associates each $E_L$ class with a countable ordinals. This suggests that $E_L$ is thin. However, the complexity of the statement that a particular $\bDelta_2^1$ equivalence relation is thin is beyond the scope of Shoenfield absoluteness. Therefore the usual argument of passing to a forcing extension satisfying $\neg \mathsf{CH}$ will not work. Morever, $E_L$ looks quite different in models that do not satisfy $\mathsf{V = L}$. Thinness will be proved more directly. 

The following fact will be useful. It implies that if $\alpha < \beta$ are admissible ordinals and a new real appears in $L_\beta$ which was not in $L_\alpha$, then $L_\alpha$ is countable from the view of $L_\beta$. 

\Begin{fact}{acceptability of J hierarcy}
If $\omega < \alpha < \beta$ are admissible ordinals and $(\cantorspace)^{L_\beta} \not\subseteq L_\alpha$, then there is an $f \in L_\beta$ such that $f : \omega \rightarrow \alpha$ is a surjection. In particular, $L_\beta \models |L_\alpha| = \aleph_0$. 
\end{fact}

\begin{proof}
This is essentially a result of Putnam. Below, a brief sketch of the proof is given using some elementary fine structure theory. (See \cite{the-fine-structure-of-the-constructible-hierarchy}, \cite{Fine-Structure}, or \cite{Introduction-to-the-Fine-Structure}.)

Note that if $\alpha$ is admissible, then $\omega \cdot \alpha = \alpha$. \cite{the-fine-structure-of-the-constructible-hierarchy} Lemma 2.15 shows that $L_\alpha = J_\alpha$, if $\alpha$ is admissible.

Now suppose $\alpha < \beta$ are admissible ordinals. Since $(\cantorspace)^{J_\beta} \not\subseteq L_\alpha$, there is some $x \in \mathscr{P}(\omega)$ such that $x \in J_\beta$ and $x \notin J_\alpha$. Then there is some $\alpha < \gamma < \beta$ and some $n \in \omega$ such that $x$ is $\Sigma_n$ definable over $J_\gamma$ but not in $J_\gamma$. \cite{the-fine-structure-of-the-constructible-hierarchy} Lemma 3.4 (i) shows that all $J_\gamma$ are $\Sigma_n$-uniformizable for all $n$. Then \cite{the-fine-structure-of-the-constructible-hierarchy} Lemma 3.1 can be applied to show that there is a $\Sigma_n$ in $J_\gamma$ surjection $f$ of $\omega$ onto $J_\gamma$. $f$ is definable in $J_\gamma$ and so $f \in J_{\gamma + 1} \subseteq J_\beta$. Since $J_\alpha \subseteq J_\gamma$, using this $f$, one can construct a surjection in $J_\beta$ from $\omega$ onto $J_\alpha$.

\end{proof}

\Begin{lemma}{greatest level less than such that real appear}
Suppose $\alpha$ is an ordinal such that there exists an $x \in \cantorspace$ with $\iota(x) = \alpha$, then there exists a greatest $\beta < \alpha$ such that there exists a $y \in \cantorspace$ with $\iota(y) = \beta$. 
\end{lemma}

\begin{proof}
Fix an $x$ such that $\iota(x) = \alpha$. If the result was not true, then there exists a sequences of reals $(x_n : n \in \omega)$ such that $\iota(x_n) < \beta$ and $\alpha = \iota(x) = \lim_{n \in \omega} \iota(x_n)$. $L_{\iota(x)} = \bigcup_{n \in \omega} L_{\iota(x_n)}$. $x \in L_{\iota(x)}$. This implies $x \in L_{\iota(x_n)}$ for some $n \in \omega$. This contradicts $\iota(x)$ being the smallest admissible ordinal $\alpha$ such that $x \in L_{\alpha}$. 
\end{proof}

\Begin{proposition}{EL is thin}
\emph{($\mathsf{V = L}$)} $E_L$ is a thin equivalence relation.
\end{proposition}

\begin{proof}
Let $T \subseteq \finBinarySequence$ be an arbitrary perfect tree. Let $\alpha = \iota(T)$. $L_\alpha$ satisfies that there are no functions from $\omega$ taking reals as images which enumerates all paths through $T$. By Lemma \ref{greatest level less than such that real appear}, let $\beta < \alpha$ be greatest such that there is a $y$ with $\iota(y) = \beta$. By fact \ref{acceptability of J hierarcy}, $L_\alpha \models |L_\beta| = \aleph_0$. However, since $L_\alpha$ satisfies no function from $\omega$ into the reals enumerate the paths through $T$, there exists $v,w \in L_\alpha$ such that in $L_\alpha$, $v$ and $w$ are paths through $T$ and $v,w \notin L_\beta$. By the choice of $\beta$, $\iota(v) = \iota(w) = \alpha$. By Proposition \ref{EL equivalent same iota}, $v \ E_L \ w$. By $\Delta_1$-absoluteness, $v, w \in [T]$. It has been shown that every perfect set has $E_L$ equivalent elements. 
\end{proof}

\Begin{remark}{drucker proof of thinness}
Motivated by Proposition \ref{EL equivalent same iota}, the above proof tries to establish the thinness of $E_L$ by studying the levels of the $L$-hierarchy. The main tool was Fact \ref{acceptability of J hierarcy}. This was proved using some fine structure theory which is somewhat technical.

Drucker has proved an equivalence relation very similar to $E_L$ is a thin equivalence relation using very simple methods from recursion theory. His method and some hyperarithmetic considerations give a far simpler proof that $E_L$ is a thin equivalence relation.
\end{remark}

\Begin{proposition}{thin equiv with all classes countable not canonicalize borel}
If $E$ is a thin equivalence relation with all classes countable, then for any $\sigma$-ideal $I$, $\{E\} \rightarrow_I \borel$ fails.
\end{proposition}

\begin{proof}
Suppose there exists some $\borel$ $I^+$ $B$ such that $E \upharpoonright B$ is $\borel$. By Silver's Dichotomy for $\coanalytic$ equivalence relations, either $E \upharpoonright B$ has countably many classes or a perfect set of pairwise $E$-inequivalent elements. The former is not possible since this would imply the $I^+$ set $B$ is a countable union of countable sets. The latter is also not possible since $E$ is thin. Contradiction.
\end{proof}

\Begin{theorem}{not EL canonicalize borel}
($\mathsf{V = L}$) For any $\sigma$-ideal $I$ on $\cantorspace$, $\{E_L\} \rightarrow_I \borel$ fails.

In particular in $L$, $\bDelta_2^{1 \ \borel} \rightarrow_I \borel$ for $\sigma$-ideal $I$ with $\bbP_I$ proper is not true.  ($\bDelta_2^{1 \ \borel}$ is the class of $\bDelta_2^1$ equivalence relation with all classes $\borel$.)
\end{theorem}

\begin{proof}
$E_L$ is thin and has all classes countable. Use Proposition \ref{thin equiv with all classes countable not canonicalize borel}.
\end{proof}

There seems to be no reason to believe that it is ever possible that the main question phrased for $\bDelta_2^1$ equivalence relations is true.

\Begin{question}{main question for pi21}
Is it consistent that for all $\sigma$-ideals $I$ such that $\bbP_I$ is proper, $\bDelta_2^{1 \ \borel} \rightarrow_I \borel$?
\end{question}

\section{Conclusion}\label{Conclusion}
This last section will put the results of this paper into perspective. Some questions will be raised and some speculations will be made. 

Large cardinal assumptions were used throughout the paper to obtain a positive answer to the main question in its various forms. In the most general case, iterability assumptions were used to get a positive answer. Iterability is a fairly strong large cardinal assumption: for example, it requires the universe to transcend $L$ in a way set forcing extensions can never do. 

However, this paper leaves open that possibility that even the most general form of this question for $\bSigma_1^{1\ \borel}$ and $\bPi_1^{1 \ \borel}$ could be provable in just $\mathsf{ZFC}$. The most interesting open question is:

\Begin{question}{negative consistency result}
It is consistent (relative large cardinals) that there is a $\sigma$-ideal $I$ on a Polish space with $\bbP_I$ proper and $E \in \bSigma_1^{1 \ \borel}$ such that $\{E\} \rightarrow_I \borel$ is false?

Same question for $\bPi_1^{1 \ \borel}$. 
\end{question}

The results of this paper provide limitations to any attempt to produce a counterexample to a positive answer to the main question. 

The results of the paper seems to suggest a universe with few and very weak large cardinals is the ideal place to consider finding a counterexample. For example, Theorem \ref{sharps of reals and positive answer main question} and Theorem \ref{sharps of reals and positive answer main question coanalytic} shows that any universe that has sharps for sets in $H_{(2^{\aleph_0})^+}$ will always give a positive answer to the main question. 

This suggest perhaps considering the question in a universe compatible with very few large cardinal, i.e. the smallest inner model of $\mathsf{ZFC}$:

\Begin{question}{status of main question in L}
What is the status of the main question in $L$? 
\end{question}

Cohen forcing ($\bbP_{I_\text{meager}}$) is perhaps the simpliest of all forcings. This paper leaves open the possibility that Cohen forcing in $L$ could be used to produce a counterexample to the main question. However, since Cohen forcing is so simple, the following is a very natural question:

\Begin{question}{cohen forcing counterexample}
Can Cohen forcing (meager ideal) be used in a counterexample to the main question? 
\end{question}

Proposition \ref{analytic equiv countable ideal arrow borel} and Proposition \ref{PIE0 equivalent silver forcing} shows that the ideal of countable sets (Sacks forcing) and the $E_0$-ideal (Prikry-Silver forcing) can never be used to produce a counterexample to the main question in the $\analytic$ case. 

One of the most common forcing extensions in descriptive set theory is the extension by the (gentle) L\'{e}vy collapse $\text{Coll}(\omega, <\kappa)$, where $\kappa$ is some inaccessible cardinal. Here there is a partial answer to Question \ref{cohen forcing counterexample}: Corollary \ref{levy collapse meager null arrow borel} shows that the meager ideal and null ideal can not be used in an extension by the L\'{e}vy collapse of an inaccessible to produced a counterexample to the main question in the $\analytic$ case. Moreover, Fact \ref{MA not CH pi12 measurable cov(I) condition} implies that these two ideals can not be used for a counterexample if $\mathsf{MA + \neg CH}$ holds.

Proposition \ref{all classes countable arrow borel} and Proposition \ref{reducible orbit equiv polish group action arrow borel} asserts that $\analytic$ equivalence relations with all classes countable or are $\borel$ reducible to orbit equivalence relations of Polish group actions can not be used to show the consistency of a negative answer. One may suspect that an unusual $\bSigma_1^{1 \ \borel}$ equivalence relation may be necessary. Thin equivalence relations include somewhat unusual objects such as $F_{\omega_1}$, $E_{\omega_1}$, and any potential counterexamples to Vaught's conjecture. However, Theorem \ref{thin analytic canonicalize ev} shows, at least in regard to the main question for $\analytic$, that thin $\analytic$ equivalence relations have the strongest form of canonicalization in the sense that one of its classes is in $I^+$.

It seems that one has reached an impasse in regard to the main question for $\bSigma_1^{1 \ \borel}$. There is a lack of interesting examples of $\bSigma_1^{1 \ \borel}$ equivalence relations which may be useful for producing a consistency result for a negative answer to the main question for $\bSigma_1^{1 \ \borel}$. 

Here is where $\bPi_1^{1 \ \borel}$ becomes much more interesting and provides a possible path forward. What appears to be promising is that $\bPi_1^1$ equivalence relations seem to be much more suspectible to set theoretic assumptions. 

One difficulty in producing the appropriate type of $\bSigma_1^1$ equivalence relation is the requirement that all classes be $\borel$. In the $\bPi_1^1$ case, one situation in which this requirement is easily solved is by considering $\bPi_1^1$ equivalence relations with all classes thin and assume $\omega_1$ is inaccessible to reals, i.e., the class $\bPi_1^{1 \ \text{thin}}$.  

Even in this case, one must still limit the universe to one in which only weak large cardinals exists: The easiest way to obtain $\omega_1$ is inaccessible to real is via a L\'{e}vy collapse. Theorem \ref{remarkable cardinal consistency pi11aleph0 arrow borel} shows that this attempt will never work if one uses a L\'{e}vy collapse extension of a remarkable cardinal. Moreover, Theorem \ref{weakly compact consistency pi11aleph0 arrow borel for c.c.c.} shows that using $\bPi_1^{1 \ \text{thin}}$ with a $\aleph_1$-c.c. forcing will never work in a L\'{e}vy collapse extension of a weakly compact cardinal. 

A closer look at the proofs of Lemma \ref{coanalytic equivalence stabilization of classes} and Lemma \ref{true sigma11 code in generic extensions} shows the follow:

\Begin{definition}{min admissible code function}
Let $E \in \bPi_1^{1 \ \borel}$. Let $r(x) = \min\{\omega_1^z : [x]_E \text{ is } \Sigma_1^1(z)\}$. 
\end{definition}

\Begin{proposition}{bound r for generic canonicalization}
Let $E \in \bPi_1^{1 \ \borel}$ and $I$ be a $\sigma$-ideal such that $\bbP_I$ is proper. Suppose for all $B \in \bbP_I$, there exists some $C \subseteq B$ with $C \in \bbP_I$ and $\sup \{r(x) : x \in C\} < \omega_1$. Then $\{E\} \rightarrow_I \borel$. 
\end{proposition}

Therefore, any counterexample to a positive answer for the main question for $\bPi_1^{1 \ \borel}$ must violate the hypothesis of this proposition. The next result gives a hypothetical condition under which this happens:

\Begin{proposition}{speculation counterexample proposition}
Suppose $\omega_1$ is inaccessible to reals. Let $I$ be a $\sigma$-ideal on a Polish space such that $\bbP_I$ is proper and whenever $g$ is $\bbP_I$-generic over $V$, $V[g] = L[g]$. Let $E \in \Pi_1^{1 \ \text{thin}}$ with the property that for all $x$, $L[x] \models [x]_E$ is uncountable thin. Then for all $C \in \bbP_I$, $\sup \{r(x) : x \in C\} = \omega_1$. 
\end{proposition}

\begin{proof}
The first claim is that $[x]_E$ can not be $\lborel(z)$ for any $z$ such that $\omega_1^z < \omega_1^{L[x]}$. (Note that $\omega_1^z$ refers to the least $z$-admissible ordinal and $\omega_1^{L[x]}$ is the least uncountable cardinal of $L[x]$.)

Suppose otherwise: $[x]_E$ is $\Sigma_1^1(z)$ and $\omega_1^z < \omega_1^{L[x]}$. As in Lemma \ref{coanalytic equivalence stabilization of classes}, define 
$$a \ E' \ b \Leftrightarrow (a \in [x]_E \wedge b \in [x]_E) \vee (a = b)$$
$E'$ is $\Sigma_1^1(z)$. $E' \subseteq E$. By the effective bounding theorem, there is some $\alpha < \omega_1^{z}$ such that $E' \subseteq E_\alpha$. Now applying Lemma \ref{pi11 increasing union borel equiv} in $L[x]$ and the fact that $\alpha < \omega_1^z < \omega_1^{L[x]}$, there exists some $\beta$ such that $\alpha < \beta < \omega_1^{L[x]}$ such that $E_\beta$ is an equivalence relation. Using the argument in Lemma \ref{coanalytic equivalence stabilization of classes}, $[x]_E = [x]_{E_\beta}$. $E_{\beta}$ is $\borel(c)$ for any $c \in \cantorspace$ such that $\text{ot}(c) = \beta$. Since $\beta < \omega_1^{L[x]}$, there exists such a $c \in L[x]$. Hence $[x]_{E_\beta}$ is $\Delta_1^1(x,c)$. 
$$V \models (\forall a)(a \ E \ x \Leftrightarrow a \ E_\beta \ x)$$
Since $x, c \in L[x]$ and this statement is $\Pi_2^1(x,c)$, by Schoenfield absoluteness
$$L[x] \models (\forall a)(a \ E \ x \Leftrightarrow y \ E_\beta \ x)$$
So $L[x] \models [x]_E$ is $\borel$. However, the assumption was that $L[x] \models [x]_E$ is uncountable thin. $\mathsf{ZFC}$ proves that no $\borel$ set can be uncountable thin. Contradiction. This proves the claim.

So now let $\alpha < \omega_1$. Let $M \prec H_\Theta$ with $\alpha \subseteq M$ and $C, \bbP_I \in M$. Note that $\omega_1^M \geq \alpha$. Let $x \in C$ be $\bbP_I$-generic over $M$. Then $M[x] \models \omega_1^{L[x]} = \omega_1^{M[x]} = \omega_1^{M} \geq \alpha$, using the fact that $\bbP_I$ has the property $V[g] = L[g]$, wherever $g$ is $\bbP_I$ generic over $V$. Certainly, the real $(\omega_1^{L[x]})^V$ is greater than or equal to $(\omega_1^{L[x]})^M \geq \alpha$. So $\omega_1^{L[x]} \geq \alpha$. By the claim above, $r(x) \geq \alpha$. Hence $\sup \{r(x) : x \in C\} = \omega_1$. 
\end{proof}

Note that if $V$ satisfies $\omega_1$ is inaccessible to reals and $V[g] = L[g]$ whenever $g$ is $\bbP_I$-generic over $V$, then ``$\omega_1$ is inaccessible to reals'' is not preserved into the extension $V[g] = L[g]$. Compare this to what happens in the $\text{Coll}(\omega, < \kappa)$ extension of $L$ when $\kappa$ is a remarkable cardinal in $L$ (see Theorem \ref{remarkable cardinal consistency pi11aleph0 arrow borel}).

Given this result, the natural questions are whether such an ideal exist and whether such an $\Pi_1^{1 \ \text{thin}}$ equivalence relation exist.

First consider the following: Suppose $\kappa \in L$ and $\kappa$ is not Mahlo. Let $G \subseteq \text{Coll}(\omega, < \kappa)$. In $L[G]$, $\omega_1^{L[G]}$ is not Mahlo  and $L[G]$ satisfies $\omega_1$ is inaccessible to reals. By \cite{Set-Theory-Exploring} Exercise 8.7, there is an $A \subseteq \omega_1$ in $L[G]$, which is reshaped, i.e., for all $\xi < \omega_1$, $L[A \cap \xi] \models |\xi| = \aleph_0$. Since $L[A] \subseteq L[G]$, $\omega_1^{L[A]} \leq \omega_1^{L[G]}$. Since $A$ is reshaped, $L[A] \models \omega_1^{L[A]} \geq \omega_1^{L[G]}$. So $\omega_1^{L[A]} = \omega_1^{L[G]}$. Since $L[G]$ satisfies $\omega_1$ is inaccessible to reals, $L[A]$ also satisfies $\omega_1$ is inaccessible to reals.

In \cite{Minimal-Coding} Section 1, it is shown that in $L[A]$ where $A$ is a reshaped subset of $\omega_1$, there is an $\aleph_1$-c.c. forcing which adds a real $g$ such that $L[A][g] = L[g]$. This forcing consists of perfect trees. By \cite{Forcing-Idealized} Corollary 2.1.5, there is a $\sigma$-ideal $I_F$ such that $\bbP_{I_F}$ is forcing equivalent to Sy-David Friedman's forcing to code subsets of $\omega_1$. In $L[A]$, $I_F$ would be a $\sigma$-ideal that satisfies the property of Proposition \ref{speculation counterexample proposition}.

It is not known whether $\sup\{r(x) : x \in C\} = \omega_1$ for all $I^+$ set $C$ is enough for a negative answer to the main question for $\bPi_1^{1 \ \text{thin}}$. It could be possible that there is a $C$ such that for all $x \in C$, $[x]_E$ is very complicated as $x$ ranges over $C$, but $C$ consists of pairwise $E$-inequivalent elements (or even $C$ is a single $E$-class). 

In $L$, Jensen's minimal nonconstructible $\Delta_3^1$ real forcing (see \cite{Definable-Sets-of-Minimal-Degree} and \cite{Set-Theory}, chapter 28) is also a forcing consisting of perfect trees. Again by \cite{Forcing-Idealized} Corollary 2.1.5, there is a $\sigma$-ideal $I_J$ such that $\bbP_{I_J}$ is forcing equivalent to Jensen's forcing. $\bbP_{I_J}$ is $\aleph_1$-c.c. by \cite{Set-Theory} Lemma 28.4. Moreover, by \cite{Set-Theory} Corollary 28.6, if $g,h$ are $\bbP_{I_J}$-generic over $L$, then $g \times h$ is $\bbP_{I_J} \times \bbP_{I_J}$ generic over $L$. Hence below any $B$ such that $B \forces_{\bbP_{I_J}} (\dot x_\text{gen})_\text{left} \ E (\dot x_\text{gen})_\text{right}$ (or $B \forces_{\bbP_{I_J}} \neg((\dot x_\text{gen})_\text{left} \ E \ (\dot x_\text{gen})_\text{right})$), if $C$ is the $I^+$ set of $\bbP_{I_J}$-generic real over $M$ in $B$ (for some $M \prec H_\Theta$), then $B$ consists of pairwise $E$-inequivalent (or pairwise $E$-equivalent) reals. But of course, this example does not satisfy all of the conditions of Proposition \ref{speculation counterexample proposition}.

It is not known whether the $\bPi_1^{1 \ \text{thin}}$ equivalence relations needed in Propopsition \ref{speculation counterexample proposition} exist.

\Begin{question}{coanalytic thin uncountable class inaccessible real question}
Let $\kappa$ be inaccessible but not Mahlo in $L$. Suppose $G \subseteq \text{Coll}(\omega, < \kappa)$ be generic over $L$. Let $A \subseteq \omega_1$ with $A \in L[G]$ be a reshaped subset of $\omega_1$. Then is there a $\Pi_1^1$ equivalence relation $E$ such that for all $x \in (\bairespace)^{L[A]}$, $L[x] \models [x]_E$ is uncountable thin?
\end{question}

This leads to an interesting related question about whether it is possible to partition $\bairespace$ in a $\Pi_1^1$ way into $\bPi_1^1$ pieces that are all uncountable thin:

\Begin{question}{partition coanalytic L}
In $L$, is there a $\Pi_1^1$ equivalence relation $E$ such that $L \models (\forall x)([x]_E \text{ is uncountable thin})$? 
\end{question}

Sy-David Friedman has communicated to the author a solution to this last question. See the appendix below for more information. 

\section{Appendix}
This appendix includes some remarks of Sy-David Friedman.

Sy-David Friedman and Tornq\"{u}ist, using some ideas of Miller and Conley, have given a solution to Question \ref{partition coanalytic L}. 

\Begin{theorem}{solution partition coanalytic L}
(Friedman, Tornq\"uist) In $L$, there exists a $\coanalytic$ equivalence relation $E$ such that $L \models (\forall x)([x]_E \text{ is uncountable thin})$. 
\end{theorem}

\begin{proof}
$E$ will be an equivalence relation on $\bbR$. Consider $\bbR$ with its usual $\bbQ$-vector space structure. By \cite{Classical-Descriptive-Set-Theory} Exercise 19.2 (i), let $C$ be a perfect $\bPi_1^0$ $\Q$-linearly independent set of reals. Let $P \subseteq C$ be an uncountable thin $\coanalytic$ subset. Let $\langle C \rangle$ and $\langle P \rangle$ denote the additive subgroups of $\R$ generated by $C$ and $P$, respectively. 

Since $C$ consists of $\Q$-linearly independent reals, each elements of $\langle C \rangle$ has a unique representation as $\Z$-linear combinations of elements of $C$. By Lusin-Novikov (countable section) uniformization, $\langle C \rangle$ is $\borel$. Also by Lusin-Novikov, there is a $\borel$ function $\Phi$ on $\R$ such that if $r \in \langle C \rangle$, then $\Phi(r)$ is a representation of $r$ as a $\Z$-linear combination of elements of $C$, and if $r \notin \langle C \rangle$, then $\Phi(r)$ is some default value.

Then $\langle P \rangle$ has the following definition: $r \in \langle P \rangle$ if and only if $r \in \langle C \rangle$ and $\Phi(r)$ consists of only elements from $P$. The latter is $\coanalytic$. Hence $\langle P \rangle$ is a coanalytic subgroup of $\R$. 

By definition, $\langle P \rangle$ is the set of $\Z$-linear combinations of elements of $P$. Since $P$ is thin, by Mansfield-Solovay, $P$ consists entirely of constructible reals. In particular, in any forcing extension $L[G]$ of $L$, $P^{L} = P^{L[G]}$. So, $\langle P \rangle^{L[G]}$ consists of $\Z$-linear combination of elements of $P^{L[G]} = P^L$. Hence, $\langle P \rangle^{L[G]} = \langle P \rangle^{L}$. If $\langle P \rangle^{L}$ had a perfect subset, then by Schoenfield's absoluteness, $\langle P \rangle^{L[G]}$ would have a perfect subset. If $G$ was generic for a forcing which makes $(2^{\aleph_0})^{L[G]} > \aleph_1^L$, then $|\langle P \rangle|^{L[G]} = (2^{\aleph_0})^{L[G]} > \aleph_1^L = |\langle P \rangle|^L$. This contradicts $\langle P \rangle^{L[G]} = \langle P \rangle$. This shows that in $L$, $\langle P \rangle$ is uncountable thin.

Let $E$ be the coset equivalence relation of $\R \slash \langle P \rangle$: $r \ E \ s \Leftrightarrow (r - s) \in \langle P \rangle$. $E$ is $\coanalytic$. For all $r$, $[r]_E$ is in bijection with $\langle P \rangle$. Hence $[r]_E$ is uncountable thin.
\end{proof}

At the time of asking Question \ref{partition coanalytic L}, there was hope that any natural constructibly coded $\bPi_1^1$ equivalence relation which witnessed a positive asnwer to Question \ref{partition coanalytic L} would also serve as a witness to a positive answer to Question \ref{coanalytic thin uncountable class inaccessible real question}.

Unfortunately, the equivalence relation $E$ of Theorem \ref{solution partition coanalytic L} does not work. The definition of $E$ has a particular constructibly coded thin $\bPi_1^1$ group built into it. $E$, as a coset relation, copies this thin uncountable (in $L$) set, throughout the reals. Now suppose $V$ is some universe such that $\omega_1^{L} < \omega_1^V$. In $V$, choose some $z \in \R$ such that $L[z] \models \omega_1^L < \omega_1$. Since $[z]_E$ is in bijection with $\langle P \rangle$ (which is in bijection with $\omega_1^L$), $L[z] \models [z]_E$ is countable. 

It seems any possible solution to Question \ref{coanalytic thin uncountable class inaccessible real question} will need to be defined without using any explicit definition of a thin $\coanalytic$ set.

\bibliographystyle{amsplain}
\bibliography{references}
\end{document}